\documentclass{article}

\usepackage{arxiv}

\usepackage[utf8]{inputenc} 
\usepackage[T1]{fontenc}    
\usepackage{hyperref}       
\usepackage{url}            

\usepackage{booktabs}       
\usepackage{amsfonts}       
\usepackage{nicefrac}       
\usepackage{microtype}      
\usepackage{graphicx}
\usepackage{doi}
\usepackage{cite}
\usepackage{amsmath,amssymb,amsfonts}
\usepackage{algorithmic}
\usepackage[ruled]{algorithm2e}
\usepackage{graphicx}
\usepackage{multirow}
\usepackage{longtable}
\usepackage[nottoc]{tocbibind}
\usepackage{array}
\usepackage{mdframed}
\usepackage{graphicx}
\usepackage{changepage}
\usepackage{soul}
\usepackage{textcomp}
\def\BibTeX{{\rm B\kern-.05em{\sc i\kern-.025em b}\kern-.08em
    T\kern-.1667em\lower.7ex\hbox{E}\kern-.125emX}}
\usepackage{ulem}
\graphicspath{ {images/} }

\newcommand{\onebf}{\boldsymbol{1}}

\newcommand{\Abf}{\boldsymbol{A}}
\newcommand{\Acal}{\mathcal{A}}

\newcommand{\cbf}{\boldsymbol{c}}

\newcommand{\dbf}{\boldsymbol{d}}

\newcommand{\Mbf}{\boldsymbol{M}}

\newcommand{\Mcal}{\mathcal{M}}

\newcommand{\Nbb}{\mathbb{N}}

\newcommand{\Rbf}{\boldsymbol{R}}

\newcommand{\Rbb}{\mathbb{R}}
\newcommand{\Rcal}{\mathcal{R}}

\newcommand{\Tcal}{\mathcal{T}}

\newcommand{\Vbf}{\boldsymbol{V}}
\newcommand{\Vcal}{\mathcal{V}}

\newcommand{\xbf}{\boldsymbol{x}}
\newcommand{\Xbf}{\boldsymbol{X}}
\newcommand{\Xcal}{\mathcal{X}}

\newcommand{\Ybf}{\boldsymbol{Y}}
\newcommand{\Ycal}{\mathcal{Y}}

\newcommand{\Zbf}{\boldsymbol{Z}}

\newcommand{\gammabf}{\boldsymbol{\gamma}}

\newcommand{\deltabf}{\boldsymbol{\delta}}

\newcommand{\diagrm}{\mathrm{diag}}

\title{$\Delta$-MILP: Deep Space Network Scheduling via Mixed-Integer Linear Programming}

\author{Thomas Claudet \\
	Jet Propulsion Laboratory\\ California Institute of Technology\\
	Pasadena, CA, 91103\\
	\texttt{thomas.claudet1@gmail.com} \\
	\And
	Ryan Alimo\thanks{Corresponding author: Ryan Alimo (e-mail: sralimo@jpl.nasa.gov).} \\
	Jet Propulsion Laboratory\\ California Institute of Technology\\
	Pasadena, CA, 91103\\
	\texttt{sralimo@jpl.nasa.gov} \\
	\And
	Edwin Goh\\
	Jet Propulsion Laboratory\\California Institute of Technology\\
	Pasadena, CA, 91103\\
	\texttt{edwin.y.goh@jpl.nasa.gov} \\
	\And 
	Mark D. Johnston\\
	Jet Propulsion Laboratory\\ California Institute of Technology\\
	Pasadena, CA, 91103\\
	\texttt{mark.d.johnston@jpl.nasa.gov} \\
	\And
	Ramtin Madani\\
	The University of Texas at Arlington\\
	TX 76015 USA\\
	\texttt{ramtin.madani@uta.edu} \\
	\And
	Brian Wilson\\
	Jet Propulsion Laboratory\\ California Institute of Technology\\
	Pasadena, CA, 91103\\
	\texttt{brian.Wilson@jpl.nasa.gov} \\
}
\date{This work was supported by the Jet Propulsion Laboratory, California Institute of Technology, under a contract with the National
Aeronautics and Space Administration (NASA). Copyrights 2021. All rights reserved.} 

\renewcommand{\shorttitle}{$\Delta$-MILP: Deep Space Network Scheduling via Mixed-Integer Linear Programming}

\begin{document}
\maketitle

\begin{abstract}
	This paper introduces $\Delta$-MILP, a powerful variant of the mixed-integer linear programming (MILP) optimization framework to solve NASA's Deep Space Network (DSN) scheduling problem. This work is an extension of our original MILP framework (DOI:10.1109/ACCESS.2021.3064928), and inherits many of its constructions and strengths, including the base MILP formulation for DSN scheduling. 
To provide more feasible schedules with respect to the DSN requirements, $\Delta$-MILP incorporates new sets of constraints including 1) splitting larger tracks into shorter segments and 2) preventing overlapping between tracks on different antennas. Additionally, $\Delta$-MILP leverages a heuristic to balance mission satisfaction and allows to prioritize certain missions in special scenarios including emergencies and landings. 
 Numerical validations demonstrate that $\Delta$-MILP now satisfies 100\% of the requested constraints and provides fair schedules amongst missions with respect to the state-of-the-art for the most oversubscribed weeks of the years 2016 and 2018.  
\end{abstract}


\section*{Nomenclature}

$\Tcal$ Set of time frames\\
$\Rcal$ Set of resources\\
$\Acal$ Set of activities\\
$\Acal^*$ Set of split activities\\
$\Acal_p$ Set of prioritized activities\\
$\Mcal$ Set of missions\\
$\Vcal$ Set of view periods\\
$\Rbf$ Mapping from resources to view periods\\
$\mathbf{A}$ Mapping from activities to view periods\\
$\mathbf{M}$ Mapping from missions to activities\\
$\mathbf{V}$ Mapping from view periods to time frames\\
$\mathbf{d}_{min}$ Minimum tracking times\\
$\mathbf{d}_{max}$ Maximum tracking times\\
$\mathbf{\delta}_{\uparrow}$ Setup times\\
$\mathbf{\delta}_{\downarrow}$ Teardown times\\
$\gammabf_{\uparrow}$ Minimum uninterrupted times once activity started\\
$\gammabf_{\downarrow}$
Minimum uninterrupted times after resource used\\
$\mathbf{x}$ Completeness of entire activities\\
$\mathbf{X}$ Allocation of view periods\\
$\mathbf{X^{\uparrow}}$ Starting time of view periods\\
$\mathbf{X^{\downarrow}}$ End time of view periods\\
$\mathbf{Y}^{\uparrow}$ Setup time interval of view periods\\
$\mathbf{Y}^{\downarrow}$ Teardown time interval of view periods\\
$\cbf_{1}$ Coefficient of priority for activities\\
$\cbf_{2}$ Coefficient of tracking time\\
$T_S$ Scheduled tracking time\\
$T_R$ Requested tracking time\\
$U_{RMS}$ Root Mean Square of the satisfaction\\
$U_{MAX}$ Most unsatisfied mission\\
$U_{AVG}$ Average satisfaction\\
$U_{PRIO}$ Satisfaction of prioritized missions\\
$d$ Metric space distance\\
$\eta_0$ Initial threshold\\
$\eta^{\uparrow}$ Threshold increment\\
$k_{max}$ Maximum number of iterations\\
$k_{time}$ Time limit of the algorithm\\

\section{Introduction}
\label{sec:introduction}
 NASA's Deep Space Network (DSN) is the primary source of communications for interplanetary spacecraft missions for most active missions around the world.  This paper considers NASA's Deep Space Network (DSN) scheduling problem, which is the problem of assigning the DSN's finite antenna resources to its users within a given time frame. DSN schedules are typically generated a
year into the future with allocations to the minute, and are performed in a semi-automated manner with manual and interactive cleanup, one week at a time \cite{Clement-2005,johnstonAutomatedSchedulingNASA2014}.
The DSN is considered as the largest and most sensitive telecommunications system in the world \cite{Taylor,NASA_management}. DSN consists of three complexes, each with one 70-meter antenna and three to four 34-meter ones \cite{Imbriale}, situated roughly 120° apart around the Earth: in Goldstone, California, Madrid, Spain, and Canberra, Australia.
NASA has estimated that deep space communications capabilities need to be expanded by a factor of 10 each decade to keep up with the rising demand for deep space communications \cite{Taylor}. 
Despite the growing demand for communications time on the DSN, its operational budget has dropped by approximately 10 million dollars per year since 2014 \cite{NASA_management}. As a result, automation plays a significant role for the efficient usage of the DSN's finite resources.
DSN is a highly utilized asset, and communications requests in recent years have frequently reached or even exceeded its full capacity. Increasing mission demands have led to frequent oversubscription --- situations in which only a subset of the requested tracking time can be feasibly allocated. For example, the Perseverance Mars rover provides higher-fidelity data compared to the previous rovers;  Perseverance includes eight cameras as well as sophisticated scientific instruments which increase the load on the DSN.

Over the past 20 years, different solutions have been proposed to automate DSN scheduling. In 2006, a genetic algorithm for DSN scheduling was proposed \cite{guillaume2007deep}, and in 2008 a generalized differential evolutionary algorithm was investigated to find a Pareto optimal solution to satisfy multiple objective functions \cite{Johnston-2008}. More recently, a deep reinforcement learning approach was proposed in which agents were trained to learn efficient strategies to generate feasible schedules based on DSN requests for a given week \cite{goh2021scheduling}.
In addition, many novel solutions and improvements have been made to increase the operational efficiency of the DSN scheduling process \cite{Bell,Chien,Johnston-2,Hackett,DSN_MILP1}.

In our previous work \cite{sabol2021deep}, the problem was formulated as a mixed-integer linear program (MILP). 
MILP has many applications such as for task graphs \cite{Roy,Roy1}, energy hub plants \cite{Jadidbonab}, flow shop scheduling \cite{Meng,Xiao2}, or job shop scheduling \cite{Meng1}. It is also extensively used for scheduling meetings \cite{Tsuruta}, to satellite operations \cite{Vazquez, Chien2,Chen,Cho,She,Wang}, or telescopes \cite{Johnston3}. In this work, we incorporated many important features of the DSN scheduling process such as spacecraft visibility/view periods, setup/teardown times before and after transmission.
Finding valid schedules that satisfy the important constraints encountered in the real world leads to a large MILP problem with tens of thousands of binary optimization variables and hundreds of millions of elements to set the problem up.

This work in \cite{sabol2021deep} presented a new, computationally efficient, and scalable formulation and algorithm to solve the DSN scheduling problem on a set of historical DSN requests. In this paper, that framework is extended by incorporating new constraints to achieve higher quality schedules in terms of feasibility and optimality. One of the primary problems with the previous approach was that the solution had less than 20\% of its tracks found to be feasible. A significant amount of scheduled tracks overlapped with one another, while others fell outside the spacecraft view periods, and other important constraints were not implemented.

In practice, the optimal schedules not only need to maximize the overall number of assigned activities but more importantly to find a balance of the satisfaction between all the missions. An activity is an input of the scheduler. It is composed of a mission name, specific setup and teardown times, and a set of view periods from which the spacecraft in visible in the sky and able to communicate with the DSN network. 
A track is an output of the scheduler. A track is a moment in time where a certain antenna will be communicating with a certain spacecraft. Tracks are then subsets of available view periods. Thus, this framework is extended to maximize the overall satisfaction across all missions while still maximizing the amount of overall tracking time, but now ensuring the respect of the realistic operational and physical constraints imposed by DSN operations.  

The main contributions of this paper can be summarized as follows:
\begin{itemize}
    \item Introducing new sets of constraints to only schedule valid tracks (tracks that satisfy all required constraints).
    \item Increasing the overall user satisfaction across missions by incorporating a dynamic objective function.
    \item Implementing a prioritization feature to maximize one or multiple mission satisfactions in case of emergencies or landings.
    \item Evaluating the new performance of our MILP solution on the most oversubscribed weeks in the DSN's operational history in 2016 and 2018.
    
\end{itemize}

The structure of this paper is as follows:  Section \S \ref{sec:jargon} introduces DSN scheduling terminology. Section \S \ref{sec:pb_formulation} derives the DSN scheduling problem formulation. Section \S \ref{sec:results} presents the results obtained by the algorithm and compares its performance with previous iterations from \cite{sabol2021deep}. Section \S\ref{sec:conclusion} finally draws conclusions on this work.
\section{Deep Space Network Scheduling Terminology} \label{sec:jargon}
This section introduces the different terms to understand how DSN scheduling works.
Many \textit{missions} today communicate with the DSN antennas. These missions are also referred as \textit{users} since they use the DSN, and antennas are also called resources. DSN produces weekly schedules, so each mission, or user, will make several communication demands to the DSN each week. These demands include a minimum and maximum tracking time during which the DSN should communicate with the spacecraft and the type of antenna that would be suited for the task. From these pieces of information, DSN first attributes for each required resource a time that will be sufficient for the antenna to prepare the communication, also called a setup time, and a time necessary after the communication also known as the teardown time. Second, knowing the position of each spacecraft at any given time in the universe with respect to the DSN antennas, DSN automatically computes the time windows during which each spacecraft is visible from each antenna. For each antenna, the window starts when the spacecraft rises from the horizon and ends when it sets. These windows are called \textit{view periods}.
The DSN then concatenates each user demand, the setup, teardown times, and the computed view periods to create what is being called an \textit{activity}. Since users make multiple requests per week, e.g., communicating with their spacecraft at least 2 hours per day, it creates as many activities.

In order to fit more activities into the schedule, the DSN allows activities to be split into smaller segments placed at different positions in the week. The rule is that activities that request more than 8 hours may be split into two 4-hour minimum chunks. The caveat to splitting activities in this fashion is that each segment now requires setup and teardown time, which is, at the end of the day, twice the amount of time that would be required if the whole track was scheduled contiguously. This thus makes split activities more costly to schedule, but can help creating space to place other tracks.

The missions operate independently of one another, and each mission would like the DSN to schedule as many of their own activities as possible. As a result, a perfect schedule from one mission's perspective will not be ideal from the point of view of all the other missions. On the other hand, the DSN wishes to efficiently utilize its resources, and has concerns such as maximizing the overall antenna usage. Thus, the goal of DSN scheduling is to satisfy the highest number of activities, while ensuring the best satisfaction balance between all the users.

\subsection{Mathematical notations}\label{sec:notations}
This section introduces the mathematical notations that are used to formulate the DSN's scheduling problem as a MILP problem. Let us denote 
$\Rbb, \Nbb, \Nbb^*$ as the sets of real numbers, natural numbers, and nonnegative integers respectively~\footnote{Calligraphic uppercase letters denote sets, mathbfface uppercase letters are used for matrices, and mathbfface lowercase for vectors.}. $|\Xcal|$ denotes the cardinality of a set $\Xcal$.
\{0,1\}$^{|\Xcal|\times|\Ycal|}$, respectively \{0,1\}$^{|\Xcal|}$ denotes the set of $|\Xcal|\times|\Ycal|$ matrices, respectively $|\Xcal|$ dimensional vectors with binary entries. $\onebf_{|\Xcal|}$, respectively $\onebf_{|\Xcal|\times|\Ycal|}$ denotes a $|\Xcal|$ dimensional vector, respectively $|\Xcal|\times|\Ycal|$ dimensional matrix of ones. $x \in \Xcal$ means that element $x$ is in the set $\Xcal$. $M_{i,j}$ represents the element on row $i$ and column $j$ of matrix $\mathbf{M}$. $\mathbf{M} \leq \mathbf {N}$ means that all the entries of $\mathbf{M}$ are less or equal to the ones of $\mathbf{N}$. $\textrm{diag}\{\Xbf\}$ is a diagonal matrix formed by the elements of vector $\Xbf$. $(.)^T$ is the transpose operator. $\mathbb{Z}/n\mathbb{Z}$ is the modular arithmetic set on the remainders of integers in the division by $n$. \{$a_1,...,a_n$\} is a list containing elements $a_1,...,a_n$. $\Xcal\subseteq \Ycal$ signifies that $\Xcal$ is a subset of $\Ycal$. $\forall x \in \Xcal$ signifies for all elements x in $\Xcal$, $\exists$ for there exist, and $\Leftrightarrow$ for equivalence.

\subsection{DSN mathematical objects definitions}

Using the notations defined in \S \ref{sec:notations}, our goal is to assign resources to spacecraft. To do so, let us define four key matrices related to DSN scheduling.

First, $\Rbf \in \{0,1\}^{|\Rcal|\times |\Vcal|}$ maps resources to view periods.  It has $|\Rcal|$ rows as well as $|\Vcal|$ columns, and for each resource $r\in \Rcal$, a 1 is assigned to the corresponding column if this resource is the one corresponding to view period $v\in \mathcal{V}$, and otherwise a 0. 

Then, $\mathbf{A}\in \{0,1\}^{|\Acal|\times |\mathcal{V}|}$ maps activities to view periods. $A_{a,v} = 1$ if and only if view period $v\in \mathcal{V}$ is found in activity $a \in \mathcal{A}$. 

Matrix $\mathbf{M}\in \{0,1\}^{|\mathcal{M}|\times |\mathcal{A}|}$ maps missions to activities. $M_{m,a} = 1$ if and only if mission $m\in \mathcal{M}$ defines activity $a \in \mathcal{A}$. 

A notion of time is finally needed to know when these view periods can be scheduled. To do that, let us finally define $\mathbf{V}\in \{0,1\}^{|\mathcal{V}|\times |\mathcal{T}|}$ that maps view periods to time frames. For example, if the time-granularity of the schedule is 1 hour, and if view period $v_i$ can only be scheduled between Monday 0:00 am and 8 am UTC, the matrix will be filled with ones the first 8 columns on the $i^{th}$ line of $\mathbf{V}$.

As stated in the previous section, a user requests a minimum and a maximum (best-case) tracking time in hours. Let us call the concatenated vector of all of them $\mathbf{d}_{min}$ and $\mathbf{d}_{max} \in \mathbb{N^{*|\Acal|}}$. 

Each antenna also has a setup and teardown time in minutes, and the concatenation of all of them across activities are: $\mathbf{\delta}_{\uparrow}$ and $\mathbf{\delta}_{\downarrow} \in \mathbb{N}^{*|\Acal|}$. To be more precise, $\mathbf{\delta_{\uparrow}}$ and $\mathbf{\delta_{\downarrow}}\in (\mathbb{Z}/15\mathbb{Z})^{|\Acal|}$ which means that the greatest common divisor of all the setups and teardowns will be equal to 15 minutes. This implies that the maximum time-granularity must be equal to 15 minutes. For example, taking it larger as 1 hour would make us lose 45 minutes at every teardown since they only last 15 minutes in reality.

Finally, $\gammabf_{\uparrow}$ and $\gammabf_{\downarrow}\in\mathbb{N}^{|\Vcal|}$ are for the minimum uninterrupted time that should be spent on each activity once it is started, and the amount of time that the resource should remain idle after the activity is interrupted.


Now, regarding the optimization variables, $\mathbf{x}\in \{0,1\}^{|\mathcal{A}|}$ tells when entire activities are completed or not throughout the planning horizon (at least the minimum requested time has been scheduled).

Also, $\mathbf{X}\in \{0,1\}^{|\mathcal{V}|\times |\mathcal{T}|}$ defines if view period $v \in \Vcal$ is allocated or not at time $t \in \Tcal$ in the schedule.

Moreover, $\mathbf{X^{\uparrow}}\in \{0,1\}^{|\mathcal{V}|\times |\mathcal{T}|}$, respectively $\mathbf{X^{\downarrow}}\in \{0,1\}^{|\mathcal{V}|\times |\mathcal{T}|}$ tells at what time $t\in \Tcal$ each view period $v \in \Vcal$ starts, respectively ends, and following \cite{sabol2021deep}:
\begin{equation}
	X^{\uparrow}_{v,t}=1 \;\Leftrightarrow\; X_{v,t-1}=0\;\wedge\;X_{v,t}=1,\label{1}
\end{equation}
\begin{equation}
X^{\downarrow}_{v,t}=1 \;\Leftrightarrow\; X_{v,t-1}=1\;\wedge\;X_{v,t}=0.\label{2}
\end{equation}

Finally, $\mathbf{Y}^{\uparrow}\in \{0,1\}^{|\mathcal{V}|\times |\mathcal{T}|}$, respectively $\mathbf{Y}^{\downarrow}\in \{0,1\}^{|\mathcal{V}|\times |\mathcal{T}|}$, tells during which interval $[t_1:t_2], \; t_1, \; t_2 \in \Tcal$ the track starts up, respectively tears down:
\begin{equation}
	Y^{\uparrow}_{v,t}=1 \;\Leftrightarrow\; \exists\;t<\tau\leq t+\delta^{\uparrow}_v: \quad X^{\uparrow}_{v,\tau}=1,\label{3}
	\end{equation}
	\begin{equation}
	Y^{\downarrow}_{v,t}=1 \;\Leftrightarrow\; \exists\;t-\delta^{\downarrow}_v<\tau\leq t: \quad X^{\downarrow}_{v,\tau}=1.\label{4}
\end{equation}


\section{Problem Statement}\label{sec:pb_formulation}
This section introduces new constraints in addition to the ones  implemented in \cite{sabol2021deep}. But all of the them are being reviewed for the purpose of completeness.

\subsection{Constraints}\label{sec:assumptions} 
The previous constraints implemented by \cite{sabol2021deep} were the following:
\begin{enumerate}
	\item Each resource can only communicate with one spacecraft at any given time. 
	\item If a track is scheduled, its tracking time must occur within a valid view period (the spacecraft must be visible in the sky above the antenna). 
	\item If a track is scheduled, it must be scheduled for at least its minimum requested tracking time, and at most its maximum requested tracking time.
	\item A track must always start with the right setup time before tracking time can occur, and end with the right teardown time. 
	\item The possibility for missions to use 2 or more antennas at the same time must be implemented (e.g. for Delta-Differential-One-Way-Ranging used for navigation measurements, or for very far away spacecraft like Voyager that need more power than what a single antenna could provide).
\end{enumerate}
These constraints correspond to equations \eqref{2b}-\eqref{2i} in the problem formulation.

In this work, new constraints are being proposed that the DSN leverages in practice.
\begin{enumerate}
    \item For the same mission, tracks must not overlap on different antennas (e.g., having a part of a first track on a certain antenna happening at the same time as another part of a second track on a different antenna):
    \begin{equation}
        \Mbf\Abf(\Ybf^{\downarrow}+\Xbf+\Ybf^{\uparrow})
\leq \onebf_{|\Mcal|\times|\Tcal|}.
    \end{equation}
Each mission $m \in \Mcal$, with $\mathbf{M}\in \{0,1\}^{|\mathcal{M}|\times |\mathcal{A}|}$ has activites $a \in \Acal$, with $\mathbf{A}\in \{0,1\}^{|\Acal|\times |\mathcal{V}|}$, who themselves have view periods $v \in \Vcal$, with $\mathbf{V}\in \{0,1\}^{|\mathcal{V}|\times |\mathcal{T}|}$ with resources $r \in \Rcal$ and available times $t \in \Tcal$ as part of them. This constraint states that, either for the setup, the tracking, or the teardown time of a track, it is forbidden for any view periods of any activities of any missions to appear more than once for a given time.
    \item Activities may be split into segments, i.e. not scheduled contiguously. 
	\begin{enumerate}
		\item If an activity is split, setup and teardown times are required for each segment. 
		\item An activity can only be split if its maximum duration is above 8 hours, and each segment must be greater than the maximum between 4 hours and half of the minimum duration. 
		\item Parts must sum up between the minimum and maximum required durations.
	\end{enumerate}
\end{enumerate}

To take into account splitting, the following algorithm based on the XOR gate in electronics is being used:
\\
\begin{algorithm}[ht]
\SetAlgoLined

 \For{each activity $a_i$}{
  \If{$d_{\mathrm{max},i} \geq 8$ hours}{
  Duplicate $a_i$ to form $a_i'$ \& $a_i''$.\\
  Append the three to the $\Acal_i^*$ set.\\
 
  $(d_{max,i'}, d_{max,i''}) \gets d_{max,i}/2$\\
  $(d_{min,i'},d_{min,i''}) \gets \max(4\; \textrm{hours}, d_{min,i}/2)$}

  Add constraints to the solver:\\
  $x_{a_i'} \leq x_{a_i''}$\\
  $x_{a_i''}  \leq x_{a_i'} $\\
  $x_{a_i}  \leq 1-x_{a_i'} $
  }
 \caption{XOR Splitting}
 \label{H} 
\end{algorithm}

Activities having a maximum duration greater or equal to 8 hours are duplicated twice. For each of them, the maximum duration is divided by 2 and the minimum duration is set to the maximum between 4 hours and half of the minimum duration. Note that the max function is being used to make the model more general in case $d_{min,i} \geq 8$ hours. Then, 3 constraints are being added to ensure that the solver either schedules the full activity or the two new duplicated activities, but nether both. One could see the electronics analogy with the Exclusive-OR (XOR - $\oplus$) gate having the truth table shown in Table \ref{tab:XOR}. 

\begin{table}[ht!]
    \centering
        \caption{\textbf{The Exclusive-OR (XOR - $\oplus$) table. P and Q are here binary variables usually used in logic. It is only desired to have the possibility to schedule either the full track $(P=1, Q=0)$, or the two duplicated tracks $(P=0, Q=1)$, but never both at the same time $(P=1, Q=1)$. Linear constraints for splitting are derived from the table and shown in Algorithm $\ref{H}$.}}
    \begin{tabular}{|c|c|c|}
    \hline
         Original track P &The two duplicated tracks Q&P$\oplus$ Q  \\
         \hline
         0&1&1\\
         \hline
         1&0&1\\
         \hline
         1&1&0\\
         \hline
    \end{tabular}\label{tab:XOR}
\end{table}
\raggedbottom
For example, if an activity requests between 8 and 10 hours, two new activities will be created with the same view periods. They will both have a minimum duration of 4 hours and a maximum duration of 5 hours. So that, if the full track is scheduled, $x_{a_i'} = x_{a_i''} = 0$ and $x_{a_i} = 1$, which verifies the inequalities stated in the algorithm, and the duration will be between the minimum and the maximum requested durations. On the other hand, if a track is split, $x_{a_i'} = x_{a_i''} = 1$ and $x_{a_i} = 0$, each duplicated track will have its minimum duration set to 4 hours and its maximum duration to 5 hours, so that the sum will have a duration between 8 and 10 hours, which is what was required.




\subsection{Standard form} \label{sec:standard_form}


The goal is to maximize the number of tracks (or scheduled view periods) as well as activities, and the objective function rewards these two quantities with the coefficients $\cbf_{1} \in\Rbb^{|\Acal|}$ and $\cbf_{2} \in\Rbb^{|\Vcal|}$.
The first coefficient $\cbf_{1}$ gives priority to the activities, when $\cbf_{2}$ rewards a better tracking time.This optimization is subject to different constraints, with \eqref{2b} - \eqref{2i} from \cite{sabol2021deep} and with the new constraints \eqref{2j} - \eqref{2m}. The complete problem formulation is thus the following:
\allowdisplaybreaks
\begin{subequations}
\begin{align}
&\underset{\begin{subarray}{c}
\Zbf
\end{subarray}}
{\text{maximize}} ~~ \cbf_1^{\top}\xbf\onebf_{|\Acal|} + \cbf_2^{\top}\Xbf\onebf_{|\Tcal|} &&& \label{2a}\\
&\mathrm{subject~to}&&&\nonumber\\
&~\Xbf\leq\Vbf &&&\label{2b}\\
&~X_{v,t}-X_{v,t-1}=X^{\uparrow}_{v,t}-X^{\downarrow}_{v,t} ~~\forall v\in\Vcal,\;\;\forall t\in\Tcal\label{2c}\\
&~\sum^{t}_{\tau=t-\gamma^{\uparrow}_v+1 }{X^{\uparrow}_{v,\tau}}\leq X_{v,t}~~~~~~~~~~\forall v\in\Vcal,\;\;\forall t\in\Tcal\label{2d}\\
&~\sum^{t}_{\tau=t-\gamma^{\downarrow}_v+1}{X^{\downarrow}_{v,\tau}} \leq 1-X_{v,t}~~~~~\forall v\in\Vcal,\;\;\forall t\in\Tcal\label{2e}\\
&~Y^{\uparrow}_{v,t}=\sum^{t+\delta^{\uparrow}_{v}}_{\phantom{-}\tau=t+1\phantom{\delta^{\downarrow}_v}}
{\!\!\!\!\!\!X^{\uparrow}_{v,\tau}}~~~~~~~~~~~~~~\forall v\in\Vcal,\;\;\forall t\in\Tcal\label{2f}\\
&~Y^{\downarrow}_{v,t}=\sum^{t}_{\tau=t+1-\delta^{\downarrow}_{v}}{\!\!\!\!\!\!X^{\downarrow}_{v,\tau}}~~~~~~~~~~~~~~\forall v\in\Vcal,\;\;\forall t\in\Tcal\label{2g}\\
&~\Rbf(\Ybf^{\downarrow}+\Xbf+\Ybf^{\uparrow})
\leq \onebf_{|\Rcal|\times|\Tcal|}&&&\label{2h}\\
&~\diagrm\{\dbf_{\min}\}\xbf\leq\Abf\Xbf\onebf_{|\Tcal|}\leq\diagrm\{\dbf_{\max}\}\xbf&&& \label{2i}\\
&~\Mbf\Abf(\Ybf^{\downarrow}+\Xbf+\Ybf^{\uparrow})
\leq \onebf_{|\Mcal|\times|\Tcal|}&&&\label{2j}\\
&~x_{a_i'} \leq x_{a_i''}~~~~~~~~~~~ \forall (a_i',a_i'') \in \Acal_i^*, ~~ \forall i\in \{1,...,|\Acal|\}\label{2k}\\
&~x_{a_i''}  \leq x_{a_i'}~~~~~~~~~~~ \forall (a_i',a_i'') \in \Acal_i^*, ~~ \forall i\in \{1,...,|\Acal|\}\label{2l}\\
&~x_{a_i}  \leq 1-x_{a_i'}~~~~~ \forall (a_i,a_i') \in \Acal_i^*, ~~ \forall i\in \{1,...,|\Acal|\}\label{2m}
\end{align}\label{eq:main}
\end{subequations}
where the  optimization variables are contained in $\Zbf$ such that:
\begin{align*}
    \Zbf := (\Xbf,\,\Xbf^{\uparrow},\,\Xbf^{\downarrow},\, 
\Ybf^{\uparrow},\,\Ybf^{\downarrow},\,\xbf).
\end{align*}
\begin{itemize}
    \item \eqref{2a} is the objective function and represents the score obtained by the total number of activities and view periods scheduled.
\item \eqref{2b} signifies that we can only schedule view periods that exist.
\item \eqref{2c} restates the definitions of \eqref{1} and \eqref{2}.
\item \eqref{2d} and \eqref{2e} enforce the minimum up and down time
constraints.
\item \eqref{2f} and \eqref{2g} define $\Ybf^{\downarrow}$ and $\Ybf^{\uparrow}$, i.e., \eqref{3} 
and \eqref{4}.
\item \eqref{2h}, analogously to \eqref{2j}, makes sure that it is impossible to assign multiple activities to a single resource at each time (including setup, communication, and teardown times). 
\item \eqref{2i} restricts the duration of each scheduled activity to be bounded by the minimum and maximum time requirements.
\end{itemize}

\subsection{User Satisfaction} \label{seq:C2}
One of the issues in \cite{sabol2021deep} which was briefly discussed was the balance between user satisfactions. It is observed that some missions were highly satisfied whereas others only had a negligible percentage of their requested hours satisfied.


Such schedules may result in mission emergencies as some missions did not get enough scheduled hours. As a result, the DSN planners will need to modify the schedules for the following week to get data from these underscheduled missions. 
In order to satisfy a certain user, the most hours has to be given to him. However, with an oversubscribed schedule, efficient trade-offs have to be made to balance the satisfaction between all users in the objective to maximize the overall one.
To compare the performances of the MILP optimization algorithms in finding an optimized solution within a specific time budget, the user satisfaction was chosen to be measured using different metrics. First, the root mean square (RMS) of the satisfaction $U_{RMS}$:
\begin{equation}U_{RMS} = \sqrt{\frac{1}{|\mathcal{M}|}\sum_{m\in \mathcal{M}}\left(\frac{T_{R_m} - T_{S_m}}{T_{R_m}}\right)^2},
\end{equation}
where $T_{S_m}$ and $T_{R_m}$ are the scheduled and requested tracking times for mission $m \in \Mcal$, and $|\mathcal{M}|$ is the number of missions.  $U_{RMS}$ measures the standard deviation of the residuals between the scheduled and the requested durations. Since it is desired to schedule as many time as requested, the residuals need to be minimized, thus the metric as well.

Then, the maximum unsatisfied fraction is defined as:
\begin{equation}
    U_{MAX} = \max_{m\in \mathcal{M}}\left(\frac{T_{R_m} - T_{S_m}}{T_{R_m}}\right).
\end{equation}
This metric measures the worst satisfaction that could be ensured for all the different missions and has to be minimized. \footnote{Minimizing $U_{MAX}$ is equivalent to maximizing the minimum satisfaction.}

And finally, the goal is to maximize the overall satisfaction ratio among missions:
\begin{equation}
    U_{AVG} = \frac{1}{|\mathcal{M}|}\sum_{m\in \mathcal{M}}\frac{T_{S_m}}{T_{R_m}}.
\end{equation}

Regarding the algorithm to improve the user satisfaction itself, the idea is to start with an objective function that does not prioritize any missions at first. A satisfaction threshold is set, which is the ratio between the requested and scheduled time of a user, and the optimization runs. After getting the results that are being saved, the satisfactions of the missions are computed. For each mission that has a satisfaction smaller than the threshold, its weights in the objective function are doubled. This is to give the mission a bigger priority as the objective function has to be maximized, and the algorithm iterates over multiples runs. If all the missions happen to have a satisfaction above the threshold after a certain iteration, the threshold is increased. If after a certain number of runs, the last threshold is not beaten, the solver stops. The last step is then to decide which result should be kept, and it will not always be the last one. It is chosen the solution that minimizes the euclidean distance $d$ in the ($U_{RMS}, U_{MAX}, 1/U_{AVG}$) space, which are the 3 chosen metrics. In case where it is asked to prioritize one or multiple missions (denoted by the $p$ index) as will be discussed in \S\ref{sub:prio}, the satisfaction metric $1/U_{PRIO}$ of these users will be added: 
\begin{equation}
    U_{PRIO} = \frac{1}{|\mathcal{M}_p|}\sum_{m_p\in \mathcal{M}_p}\frac{T_{S_{m_p}}}{T_{R_{m_p}}}.
\end{equation}
The general-case distance $d$ is thus the following:
\begin{equation}
d = \sqrt{U_{RMS}^2 + U_{MAX}^2 + \frac{1}{U_{AVG}^2} + \frac{1}{U_{PRIO}^2}}.
\end{equation}

Note that the two first have to be minimized, when the two last have to be maximized. This is why the inverse is taken for them. Algorithm \ref{H2} summarizes this.

\begin{algorithm}[h]
\SetAlgoLined
$k \gets 0$\\
$\eta \gets \eta_0$ \\
\While{$k < k_{max}$}{
    Solve problem $\eqref{2a}-\eqref{2m}$ during $k_{time}$\\
    Save the solution.\\
    \For{Each mission in the results}{
        Compute its satisfaction (ratio of scheduled hours by requested).\\
        \If{satisfaction < $\eta$}
            {Double its weights in the objective function.}
    }
    \While{min(satisfactions) $\geq \eta$}
        {
            $\eta \gets \eta + \eta^+$\\
            $k \gets 0$
        }
    \If{Current solution = Previous solution}{
        $k_{time} \gets  2.k_{time}$\\
        $k \gets 0$
    }
    $k \gets k + 1$
}  
\For{Each solution}{
    Compute $U_{RMS}, U_{MAX}, 1/U_{AVG}, 1/U_{PRIO}$. \\
    Compute distance d and keep the solution that minimizes it.
}
\caption{Dynamic objective function}
\label{H2}
\end{algorithm}

\section{Numerical Results}
\label{sec:results}

\subsection{Input Datasets} \label{sec:inputdata}
A set of User Loading Profiles (ULPs) was used for Week 44 in 2016 (W44 2016) as well as weeks 10, 20, 30, 40, and 50 of 2018, all of which being oversubscribed. The algorithm thus has to make efficient choices to increase the chosen metrics. 

The input data is the set of all the activities for these given weeks. Table \ref{tab:missiondata2} summarizes the dataset that is used in this work. 
\begin{table}[ht!]
\caption{\textbf{Summary of the input data used for DSN scheduling problems. \# Resources represents the number of available antennas, Requested $n_a$ is the total number of requested activities during that week, $T_R$ is the total tracking time in hours requested from all missions, and \# Missions is the number of missions asking for hours during that specific week.}}
\centering
\scalebox{0.85}{
\begin{tabular}{|c|c|c|c|c|}
    \hline
     \textbf{Dataset} & {\# Resources} & Requested $n_a$ & $T_R$ (hrs) & \# Missions   \\
\hline     
\textbf{W44 2016} & 14 & 284 & 1418 & 29  \\
\hline
\textbf{W10 2018} & 12 & 257 & 	1192 & 	30  \\
\hline
\textbf{W20 2018}  & 12 & 294 & 	1406 & 	33  \\
\hline
\textbf{W30 2018}  & 12 & 293& 	1464& 	32  \\
\hline
\textbf{W40 2018}  & 12 & 333 & 	1737 & 	33  \\
\hline
\textbf{W50 2018}  & 12 & 275 & 	1292 & 	29  \\
\hline
\end{tabular}}
 \label{tab:missiondata2}
\end{table}

This section evaluates the new framework on real-world datasets from NASA's Deep Space Network historical schedules. Numerical evaluations were performed on these 6 oversubscribed weeks and results were compared with the state-of-the-art.

\subsection{Experimenting on week 44, 2016}\label{sec:w442016}
Week 44 in 2016 (W44 2016) is an oversubscribed week where several algorithms tried to find the optimal schedule for, including  \cite{goh2021scheduling} with Deep Reinforcement Learning (DeepRL) and \cite{sabol2021deep} with Mixed-Integer Linear Programming (MILP).  As shown in Table \ref{tab:missiondata2}, W44 2016 consists of 14 resources, 284 activities, 1418 hours of requested tracking time, and 29 different missions. Also, Table \ref{tab:missiondata_week44}, provides a summary of the parameters associated with the 284 activities in this data set. Each week in our dataset is discretized into a set of 15-minute time intervals.  

\begin{table}[ht!] 
\caption{\textbf{Average parameter values organized by missions for the week 44 in 2016 DSN data set.
$T_R$ is the total number of requested hours for each user. $n_a$ is the number of activities for each user. $d_{min}$ and $d_{max}$ are the average minimum and maximum durations requested by the user in hours. $\delta^{\uparrow}$ and $\delta^{\downarrow}$ are the average setup and teardown times of all the requests of the user in minutes.}}
\centering
\scalebox{0.85}{
\begin{tabular}{|c |c |c|c|c|c|c|}
    \hline
     \textbf{Mission} & \textbf{$T_R$} (hrs) & \textbf{$n_a$} & $\dbf_{\min}$(hrs) & $\dbf_{\max}$(hrs) & $\deltabf_{\uparrow}$(mins) & $\deltabf_{\downarrow}$(mins)  \\
\hline     
\textbf{ACE} & 27.5 & 10 & 2.75 & 2.75 & 60 & 15 \\
\hline
\textbf{CAS} & 72.0  & 8 & 7.2 &	9 &	60 & 15 \\
\hline
\textbf{CHDR} & 21 & 21 & 1 &	1 &	60 & 15 \\
\hline
\textbf{DAWN}  & 168 & 21 & 6.4 & 8 & 60 & 15 \\
\hline
\textbf{DSCO} & 1 & 1 &	1 &	1 & 60 & 15 \\
\hline
\textbf{GTL} & 23.1	& 21 &	1.1 & 1.1 &	30 & 15 \\
\hline
\textbf{HYB2}& 6 &	2 &	3 & 3 & 60 & 15 \\
\hline
\textbf{JNO} & 128	& 18 & 5.95 & 7.11 & 63.33 & 16.67 \\
\hline
\textbf{KEPL} & 12 & 2 & 6 & 6 & 60 & 15 \\ 
\hline
\textbf{LRO} & 13	& 13 &	1 &	1 &	60 & 15 \\
\hline
\textbf{M01O} & 113 &	15 & 6.04 & 7.53 & 60 & 15 \\
\hline
\textbf{MER1} & 17.5 &	7 &	2.5 & 2.5 &	45 & 15 \\
\hline
\textbf{MEX} & 56	 & 7 & 6.4 & 8 & 60 & 15 \\
\hline
\textbf{MMS1} & 51.2 &	13 & 3.94 &	3.94 &	60 & 15 \\
\hline
\textbf{MOM} & 25.5	 & 8 & 3.19 & 3.19 & 46.88 & 15 \\
\hline
\textbf{MRO} & 112	 & 14 &	6.4 & 8 & 60 & 15 \\
\hline
\textbf{MSL} & 42	& 7	& 6 & 6 & 60 & 15 \\
\hline
\textbf{MVN} & 72	& 9 & 6.4 &	8 &	60 & 15 \\
\hline
\textbf{NHPC} & 70 &	7 &	8 &	10 & 60 & 15 \\
\hline
\textbf{ORX} & 45	& 10 &	4.5 & 4.5 & 60 & 15 \\
\hline
\textbf{SOHO} & 84 &	21 & 4 & 4 & 60 & 15 \\
\hline
\textbf{STA} & 28	& 7 & 4 & 4 & 60 & 15 \\
\hline
\textbf{STB} & 15.5	& 3 & 5 & 5.17 & 45 & 15 \\
\hline
\textbf{STF} & 20	& 5	& 4 & 4 & 60 & 15 \\
\hline
\textbf{THB} & 14	& 4	& 3.5 &	3.5 & 60 & 15 \\
\hline
\textbf{THC} & 14	& 4	& 3.5 &	3.5 & 60 & 15 \\
\hline
\textbf{VGR1} & 65 &	8 &	6.67 &	8.13 & 46.88 & 15 \\
\hline
\textbf{VGR2} & 84.5 &	11 & 6.41 &	7.68 &	35.45 &	15 \\
\hline
\textbf{WIND} & 17.5 &	7 &	2.5 & 2.5 &	60 & 15 \\
\hline
\end{tabular}}
\label{tab:missiondata_week44}
\end{table}
\nobreak

One of the challenges of the DSN's scheduling problem is the nonuniform requested time distribution between missions. Table \ref{tab:missiondata_week44} illustrates that for week 44, 2016.  For instance, DAWN\footnote{DAWN is a mission studying Vesta and Ceres, two protoplanets in the asteroid belt.} requested 168 hours of tracking time while DSCO~\footnote{DSCO is the Deep Space Climate Observatory} only asked for 1 hour. As a result, a single hour of tracking time for DSCO is a lot more important than 1 hour for DAWN. 

\begin{figure}[ht!]
    \centering
  \includegraphics[width=1\textwidth]{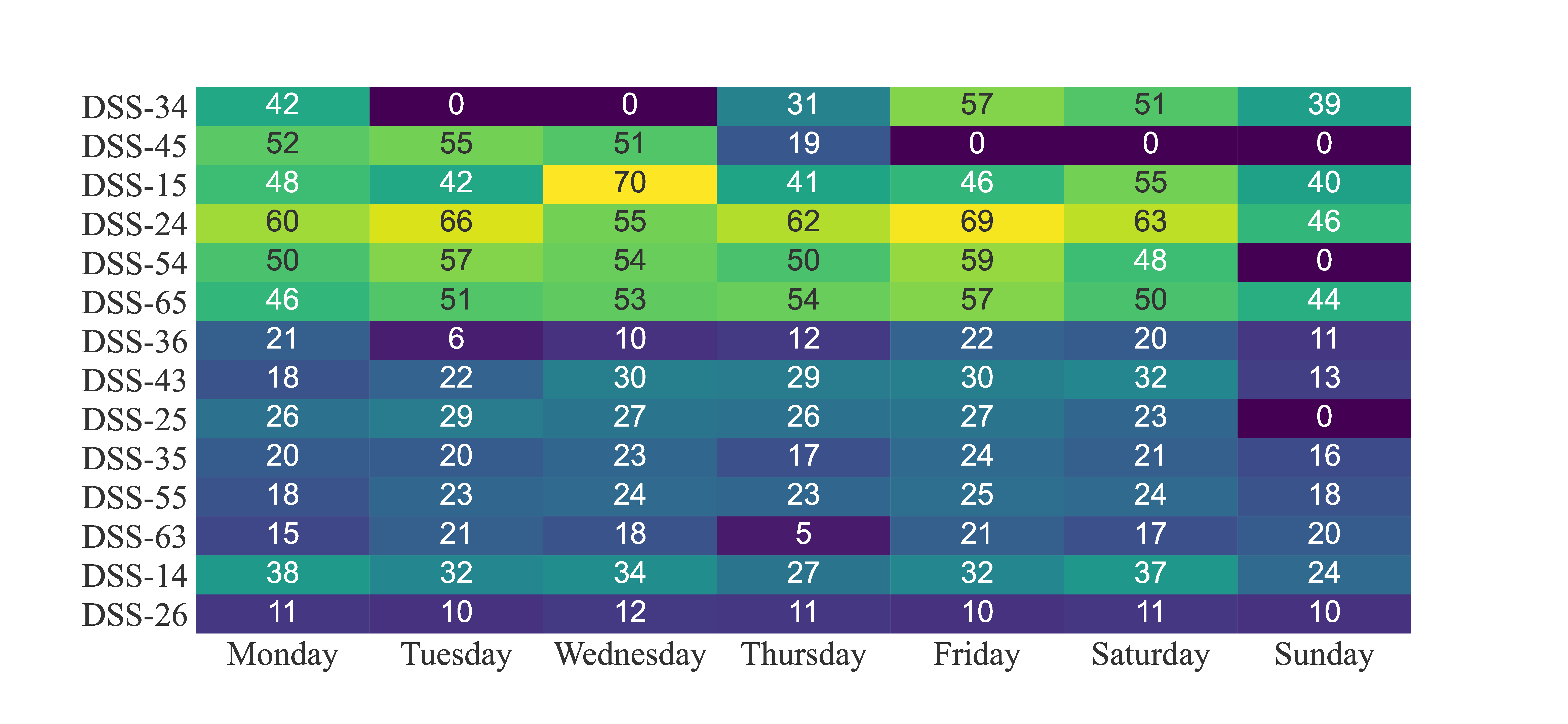}
    \caption{\textbf{Availability of different antennas for the different days of the week. For each cell, the number of times the related antenna could be requested on a given day is indicated. It is observed that there is a non-homogeneity between antennas. Some of them were under maintenance or got completely shut down such as with DSS-34 and DSS-45. Moreover, since DSN antennas do not have the same diameter and are not placed at the same location, there are some disparities between antennas of a same complex on a given day. This enforces the difficulty of balancing mission satisfactions in this setup.}}
    \label{fig:antenna_competition_W44_2016}
\end{figure}

Using $\Delta$-MILP, for week 44 in 2016, Fig. \ref{fig:antenna_competition_W44_2016} illustrates the availability of different antennas for the different days of the week. For each cell, the number of times the related antenna could be requested on a given day is indicated. It is observed that there is a non-homogeneity between antennas. Some of them were under maintenance or got completely shut down such as with DSS-34 and DSS-45. Moreover, since DSN antennas do not have the same diameter and are not placed at the same location, there are some disparities between antennas of a same complex on a given day. This enforces the difficulty of balancing mission satisfactions in this setup.

\begin{figure}[ht!]

            \includegraphics[width=1\columnwidth]{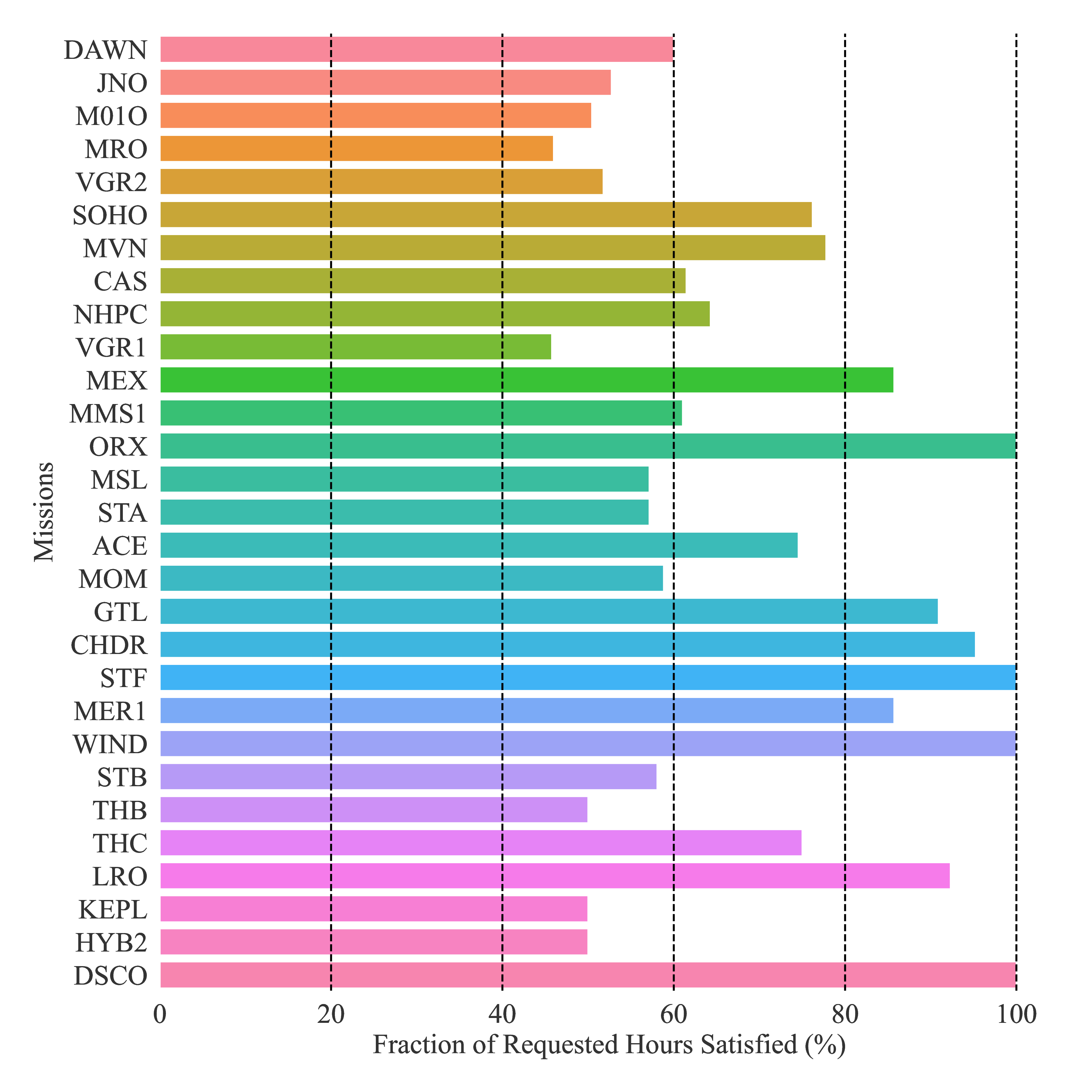}
    \caption{\textbf{Satisfaction ratios for W44 2016, the most oversubscribed week of the year, using $\Delta$-MILP. 100\% valid tracks were scheduled. 901 hours (63.5\%) as well as 208 requests (73.2\%) were satisfied. The average satisfaction ratio $U_{AVG}$ is 69.9\%, the $U_{RMS}$ is 0.35 and the minimum satisfaction is 45.8\% $(U_{MAX} = 54.2\%)$ which strongly outperforms any other implementations developed so far.}}
      \label{fig:results_W44_2016_MILP3}
        \end{figure}

The following bar chart on Fig. \ref{fig:results_W44_2016_MILP3} gives the satisfaction ratios for W44 2016, the most oversubscribed week of the year, using $\Delta$-MILP. 100\% valid tracks were scheduled. 901 hours (63.5\%) as well as 208 requests (73.2\%) were satisfied. The average satisfaction ratio $U_{AVG}$ is 69.9\%, the $U_{RMS}$ is 0.35 and the minimum satisfaction is 45.8\% $(U_{MAX} = 54.2\%)$ which strongly outperforms any other implementations developed so far.
All MILP formulations are solved using GUROBI 9.1.2 through MATLAB R2020b on a laptop computer with Intel Core i7 4.1 GHz CPU, and 16 GB of RAM. 
Regarding the initialization of Algorithm \ref{H2}, the objective function was initialized by giving an equal weight to every completed activity and assigned view period ($c_1 \gets \mathbf{1}_{|\Acal|}$, $c_2 \gets \mathbf{1}_{|\Vcal|}$ from (\ref{2a})). The maximum number of iterations $k_{max}$ was assigned to 10 to have enough runs to efficiently balance the satisfaction. The time limit per iteration $k_{time}$ was set to 30 minutes as the solver generally found and for some weeks converged to a solution after that time. If it cannot, following Algorithm \ref{H2}, the time limit doubles. The initial threshold $\eta_0$ was set to 15\% and its increment $\eta^+$ to 5\%.

Results for $\Delta$-MILP show that the satisfaction is well balanced between every mission with a minimum of 45.8\% provided by Algorithm \ref{H2}. Only valid tracks were scheduled and more than 63\% of hours and  73\% of requests were satisfied. These numbers are significantly higher than the state-of-the-art as shown in Table \ref{tab:comparision_weeks}.
        
\begin{figure*}[ht!]
    
  \includegraphics[width=\textwidth]{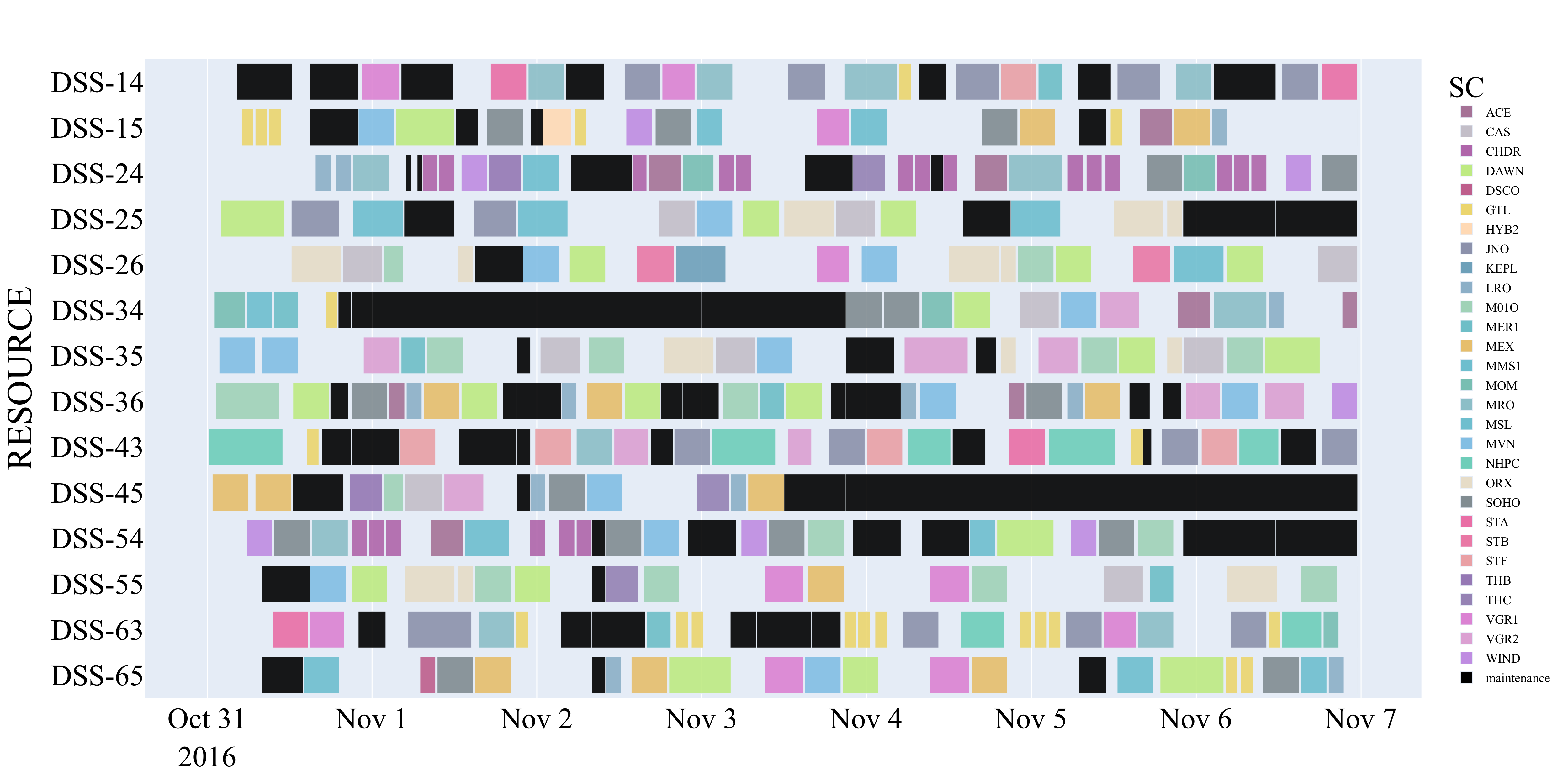}
    \caption{\textbf{Final schedule using $\Delta$-MILP. All missions have a different color, and a technician hovering them with a mouse would see more precisions such as the mission's name, the exact moment of the setup, the communication, and the tear-down time, but also the track-id and the antennas involved. There are still some time frames during which no tracks were placed. This is either because no view periods were available at that time to be scheduled, or since it did not produce a better solution. SC means Spacecraft.}}
    \label{fig:schedule_MILP3_W44_2016}
\end{figure*}

On Fig. \ref{fig:schedule_MILP3_W44_2016}, the final schedule using $\Delta$-MILP is shown. All missions have a different color, and a technician hovering them with a mouse would see more precisions such as the mission's name, the exact moment of the setup, the communication, and the tear-down time, but also the track-id and the antennas involved. There are still some time frames during which no tracks were placed. This is either because no view periods were available at that time to be scheduled, or since it did not produce a better solution.

\begin{figure}[ht!]
  \hspace*{-0cm}\includegraphics[width=1\textwidth]{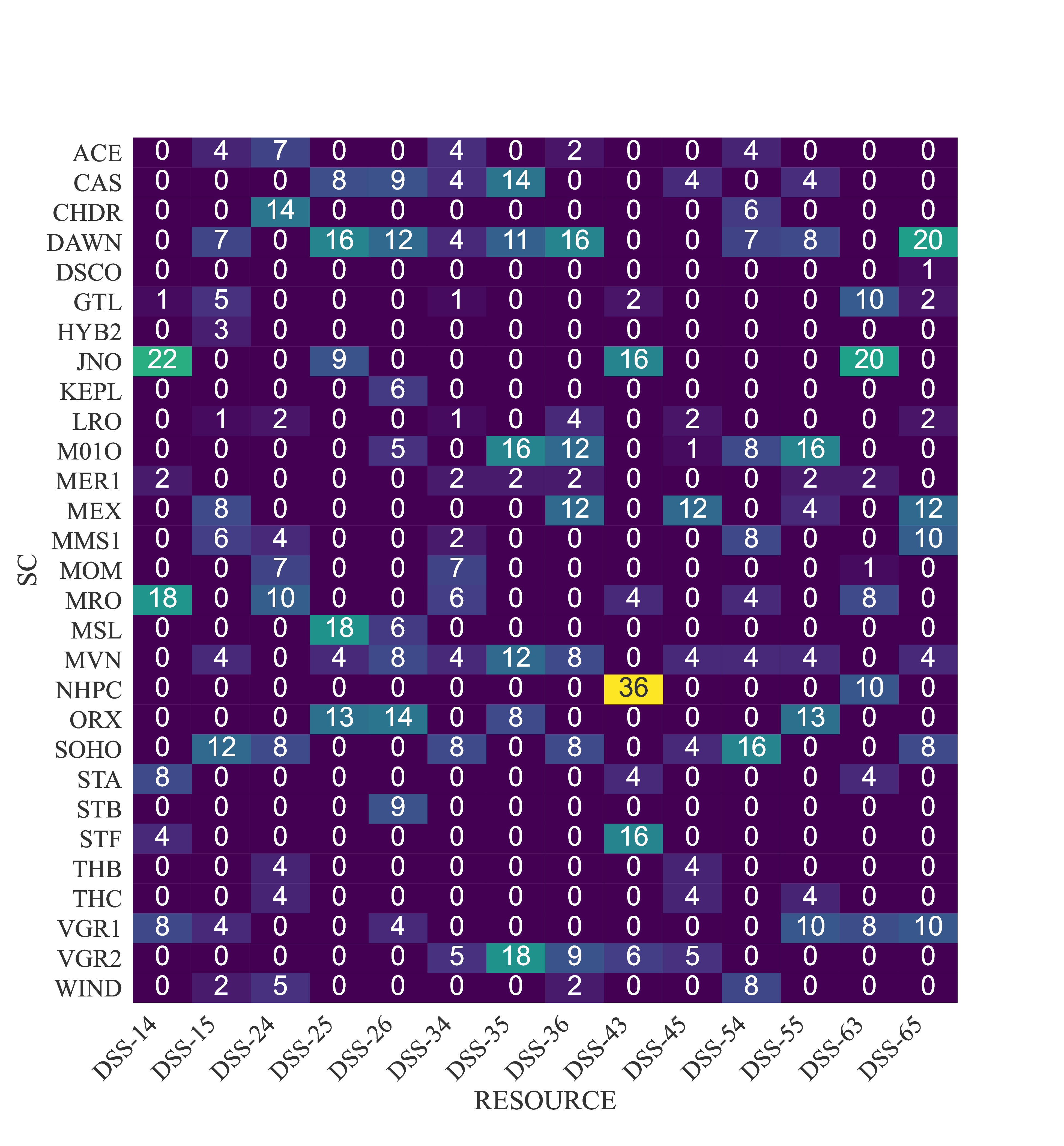}
    \caption{\textbf{Heatmap of the scheduled hours for each mission and on each antenna is displayed. The cell for NHPC on DSS-43 is very bright since the mission only required 70-m antennas, and DSS-43 was the only one that had view periods with the spacecraft during that week.}}
    \label{fig:heatmap_scheduled_durations_antennas_MILP_W44_2016}
\end{figure}

Then, on Fig. \ref{fig:heatmap_scheduled_durations_antennas_MILP_W44_2016}, a heatmap of the scheduled hours for each mission and on each antenna is displayed. The cell for NHPC \footnote{New Horizon is a mission to study the dwarf planet Pluto.} on DSS-43 is very bright since the mission only required 70-m antennas, and DSS-43 was the only one that had view periods with the spacecraft during that week. 

\begin{figure}[ht!]
  \hspace*{-0cm}\includegraphics[width=1\textwidth]{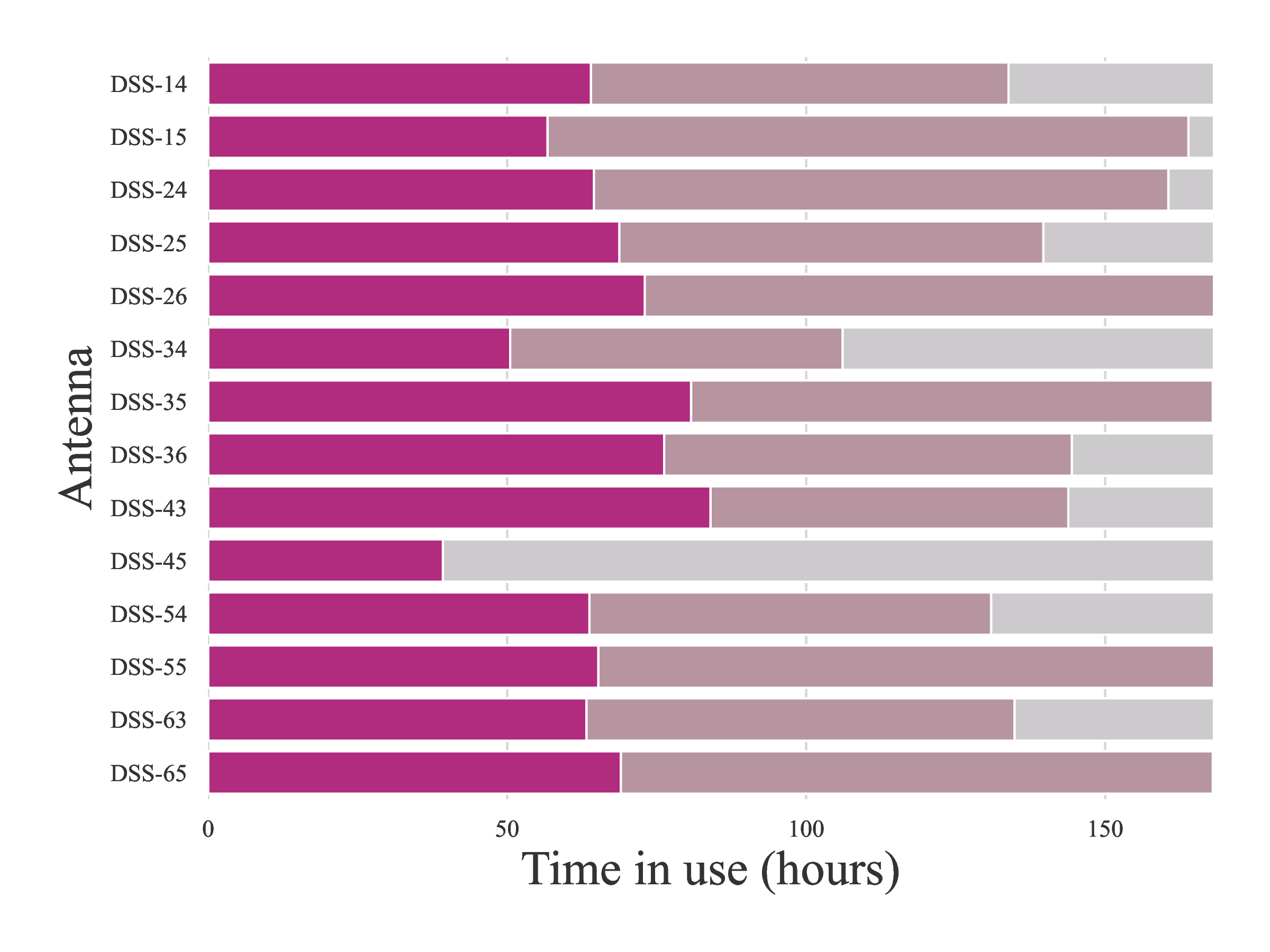}
    \caption{\textbf{Antenna usage for W44 2016. Communication time is in red, available time in light red, and maintenance or shutdown time in grey. DSS-43 scheduled all of its available hours. This bar chart shows a balance between all the hours of the different antennas, which is one of the objectives from the DSN perspective.}}
    \label{fig:antennas_use_MILP3_W44_2016}
\end{figure}

Finally, on Fig. \ref{fig:antennas_use_MILP3_W44_2016}, the antenna usage for W44 2016 is shown. Communication time is in red, available time in light red, and maintenance or shutdown time in grey. DSS-43 scheduled all of its available hours. This bar chart shows a balance between all the hours of the different antennas, which is one of the objectives from the DSN perspective.


\subsection{More comparisons with weeks 10, 20, 30, 40, and 50 in 2018}\label{sec:w10203040502016}

In addition to week 44 in 2016, the real world DSN's dataset was also studied to propose new comparisons for weeks 10, 20, 30, 40, and 50 in 2018. Results are shown in Table \ref{tab:comparision_weeks}. 4 algorithms were compared. The first two were the ones introduced in \cite{sabol2021deep} using an unweighted objective function as well as randomizing the weights at every iteration. In \cite{sabol2021deep}, they were uniformly drawn from the interval $[1,5]$ and $[0,0.01]$, respectively. Then, there is $\Delta$-MILP, and the results after its first iteration (30 minutes of run time), $\Delta$-MILP(0).

For every week, it is observed that $\Delta$-MILP outperforms by more by 5 times the original MILP \cite{sabol2021deep} in terms of hours satisfied, overall satisfied time fraction, number of satisfied requests, as well as overall satisfied request fraction, average satisfied ratios, RMS and maximum unsatisfied time fraction. Algorithms provided by \cite{sabol2021deep} did not respect all the constraints, that is why the goal of this work also was to extend the constraints set. It should be noted that no more constraints were added to be assessed between \cite{sabol2021deep} and $\Delta$-MILP.

It can also be seen that $\Delta$-MILP(0) sometimes places more hours or satisfy more requests than $\Delta$-MILP. However, when looking at the $U_{MAX}$ metric, the satisfaction of $\Delta$-MILP outperforms $\Delta$-MILP(0), showing the usefulness of Algorithm \ref{H2} to increase the minimum satisfaction using this dynamic objective function and iterative threshold approach.

\begin{table*}[ht!]
    \centering
    \caption{\textbf{Comparison of scheduling results using previously implemented algorithms and new ones. A valid track is a track that satisfies all the constraints defined in  \S \ref{sec:standard_form}}}
    \centering
    \renewcommand{\thetable}{\Roman{table}}
    \begin{tabular}{ |c|c|c|c|c|c|c|c|  }
 \hline
 Metric & Alg.  & W44\_16 & W10\_18 & W20\_18 & W30\_18& W40\_18& W50\_18\\
 \hline
 
 \multirow{4}{*}{Valid tracks (\%)} 
 &   Unweighted \cite{sabol2021deep}   & 16.5& 19.0  &  14.2 & 11.9   & 13.5 & 14.3   \\ \cline{2-8}
  & Randomized \cite{sabol2021deep}&22.1 &20.2 &22.0 &40.0 & 32.9& 25.1\\ \cline{2-8}
 & $\Delta$-MILP(0) &100 &100&100 & 100&100 &100 \\\cline{2-8}
  &  $\Delta$-MILP& 100 &100   & 100  & 100  & 100 & 100      \\
    \hline

\multirow{4}{*}{Hours satisfied} 
 &   Unweighted \cite{sabol2021deep}  & 171 & 152   &  136 &  135  & 114 & 127  \\ \cline{2-8}
  & Randomized \cite{sabol2021deep}& 207&159 &156 & 127& 70& 134\\ \cline{2-8}
    & $\Delta$-MILP(0) & 900& 855&990 & 990& 885& 821\\ \cline{2-8}
  &  $\Delta$-MILP &901 &  822  &  1059 & 983  & 949 &   816    \\
    \hline
      
\multirow{4}*{\begin{tabular}{l}
                   Overall satisfied \\time fraction (\%)\end{tabular}}
 &   Unweighted \cite{sabol2021deep}  &11.9  & 12.2  &  9.7 &  9.2  & 6.5 & 9.8 \\ \cline{2-8}
    & Randomized \cite{sabol2021deep}&14.6 &13.3 &11.1 & 8.7&4.1 &10.4 \\ \cline{2-8}
    & $\Delta$-MILP(0)& 63.5 & 71.8& 70.4&67.6 & 51.0&63.5 \\\cline{2-8}
  &  $\Delta$-MILP &63.5 & 69.1  & 75.4&67.1  & 54.6 & 63.1 \\
    \hline
\multirow{4}{*}{\# of satisfied requests} 
 &   Unweighted \cite{sabol2021deep}  & 51 & 47   &  38& 34 & 43  &41\\ \cline{2-8}
    & Randomized \cite{sabol2021deep}&62 & 40&46 & 31&24& 45\\ \cline{2-8}
    & $\Delta$-MILP(0)& 216 &212& 239& 232& 223&212 \\\cline{2-8}
  &  $\Delta$-MILP&  208& 203  & 249 & 231 & 223 &197  \\
    \hline
\multirow{4}{*}{\begin{tabular}{l}
                   \;Overall satisfied \\ request fraction (\%)\end{tabular}} 
 &   Unweighted \cite{sabol2021deep}  & 17.8 &17.5   & 12.9&11.5 &  12.6 & 14.7 \\ \cline{2-8}
    & Randomized \cite{sabol2021deep}&21.8 &15.6 & 15.6& 10.6&7.2 &16.4 \\ \cline{2-8}
    & $\Delta$-MILP(0)& 76.1 &82.5 &81.3 &79.2 & 67.0&77.1 \\\cline{2-8}
  &  $\Delta$-MILP &73.2 &  79.0  & 84.7& 78.8 &  67.0 & 71.6  \\
    \hline
\multirow{4}{*}{$U_{AVG}$ } 
 &   Unweighted \cite{sabol2021deep}  & 21.2 & 22.7   & 17.9& 9.4 &  17.1  & 17.4 \\ \cline{2-8}
    & Randomized \cite{sabol2021deep}& 22.9&20.1 &18.0 &8.9 &5.5 & 14.5\\ \cline{2-8}
    & $\Delta$-MILP(0)& 76.3 & 85.4& 86.7&82.5 & 70.3&81.3 \\\cline{2-8}
  &  $\Delta$-MILP &69.9 &  81.5  & 88.8& 81.4 & 70.8 & 73.8 \\
    \hline
\multirow{4}{*}{$U_{RMS}$ } 
 &   Unweighted \cite{sabol2021deep}  &0.82  &0.82   & 0.86& 0.92 &0.88  & 0.85 \\ \cline{2-8}
    & Randomized \cite{sabol2021deep}& 0.80& 0.84& 0.85 & 0.92& 0.95&0.87 \\ \cline{2-8}
    & $\Delta$-MILP(0)& 0.36 & 0.27& 0.25& 0.28& 0.44& 0.30\\\cline{2-8}
  &  $\Delta$-MILP & 0.35 & 0.26 & 0.21& 0.29&   0.40&0.35   \\ 
  \hline
  \multirow{4}{*}{$U_{MAX}$} 
 &   Unweighted \cite{sabol2021deep} & 100  & 100  &100 &100  &100  &100 \\ \cline{2-8}
    & Randomized \cite{sabol2021deep}&100 & 100&100 & 100& 100&100 \\ \cline{2-8}
    & $\Delta$-MILP(0)& 86.4 & 91.7& 72.7&79.2 & 100& 72.8\\\cline{2-8}
  &  $\Delta$-MILP &54.2& 47.9 & 64.1&64.3 &100  &60.0  \\
    \hline
    \multirow{4}{*}{Run time (hours)} 
 &   Unweighted \cite{sabol2021deep} & 15 &  15 &15 &  15& 15 & 15\\ \cline{2-8}
    & Randomized \cite{sabol2021deep}& 15& 15& 15&15 & 15& 15\\ \cline{2-8}
    & $\Delta$-MILP(0)& 0.5 & 0.5& 0.5&0.5 & 0.5&0.5\\\cline{2-8}
  &  $\Delta$-MILP &7.5&18 &10.5& 13.5&22.5 & 7.5 \\
  \hline
 \end{tabular}

    \label{tab:comparision_weeks}
\end{table*}

\begin{figure}[ht!]
  \hspace*{-0cm}\includegraphics[width=1\textwidth]{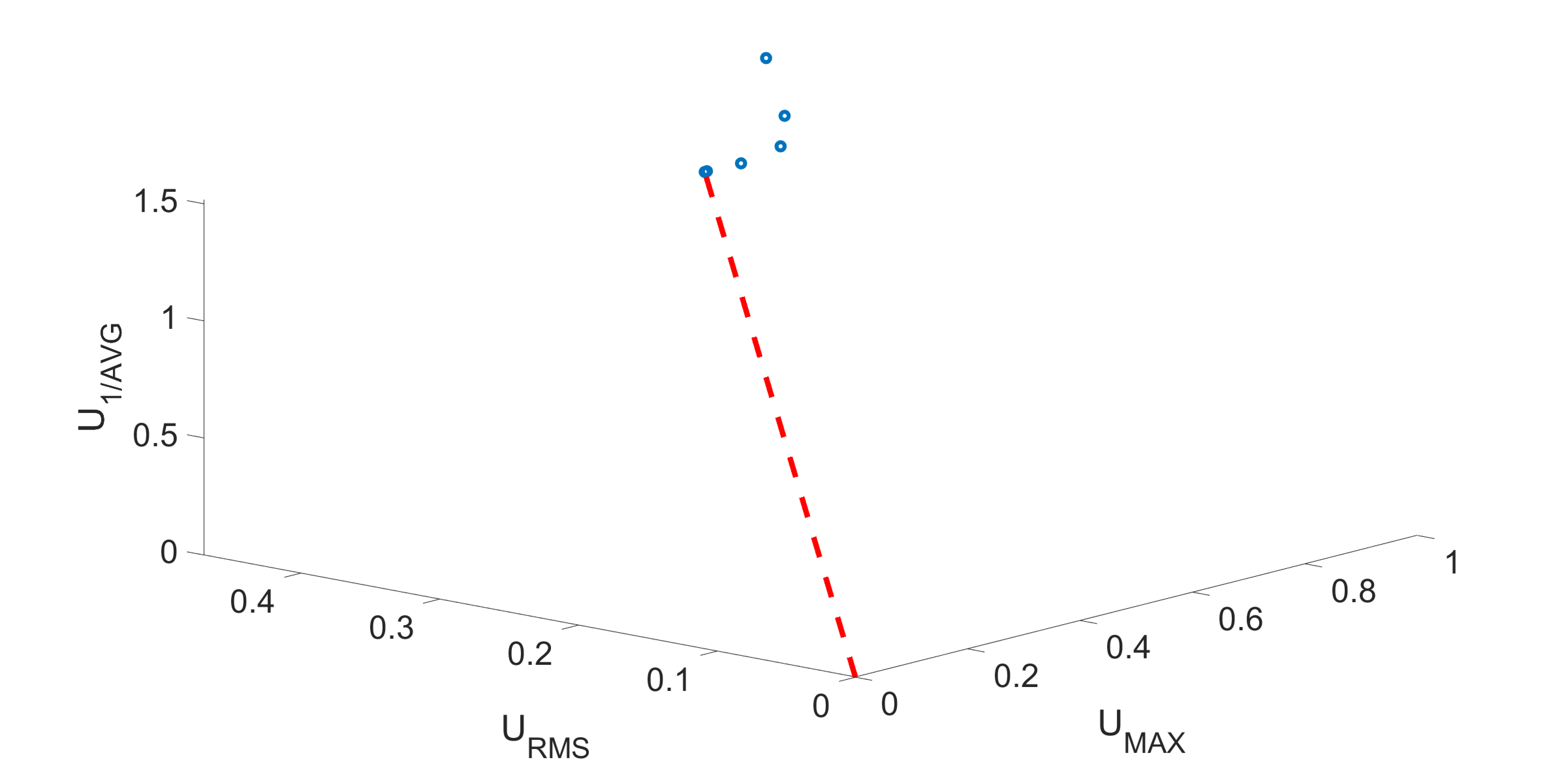}
    \caption{\textbf{The 3 metrics $(U_{RMS}, \; U_{MAX}$, and $1/U_{AVG})$ computed from the iterations of $\Delta$-MILP from Algorithm \ref{H2} for W40 2018 which is the hardest of this year. The blue dots are the iterative solutions and the red broken line indicates the minimum distance, and thus the chosen schedule in this metric space. The schedules gradually improve with respect to them.}}
    \label{fig:iterations}
\end{figure}

Fig. \ref{fig:iterations} shows the 3 metrics $(U_{RMS}, \; U_{MAX}$, and $1/U_{AVG})$ computed from the iterations of $\Delta$-MILP from Algorithm \ref{H2} for W40 2018 which is the hardest of this year. The blue dots are the iterative solutions and the red broken line indicates the minimum distance, and thus the chosen schedule in this metric space. The schedules gradually improve with respect to them.
\subsection{Prioritizing a mission} \label{sub:prio}
A last feature was finally implemented to give the possibility to prioritize one or multiple missions. Prioritization is a key factor during the DSN operations. Examples of which would be the landing of a new spacecraft or an emergency where there is a chance to lose the spacecraft (e.g.,  the Spirit rover from its Martian sand trap \cite{Spirit}). In such cases, our algorithm for DSN scheduling has to create a trade-off to be able to schedule as many hours as possible for these specific missions while not decreasing the satisfaction of the other ones too much.

When incorporating the priority requests for the DSN, the overall implementation will remain as it was described in \S \ref{sec:standard_form} for Algorithm \ref{H2}. However, instead of starting with weights of ones within the objective function as given in \eqref{2a}, it is considered a nonuniform vector that embeds the priority for each mission. The larger the multiplier is, the more important the mission will become. Here, after discussing with the DSN operators and considering their experience in this area,  it was considered a multiplier of 5 for the missions with a higher priority, i.e. $c_1(i_p) \gets 5, \; \forall i_p \in \Acal_p$, the prioritized activities set. 

\begin{figure}[ht!]
  \hspace*{-0cm}\includegraphics[width=1\textwidth]{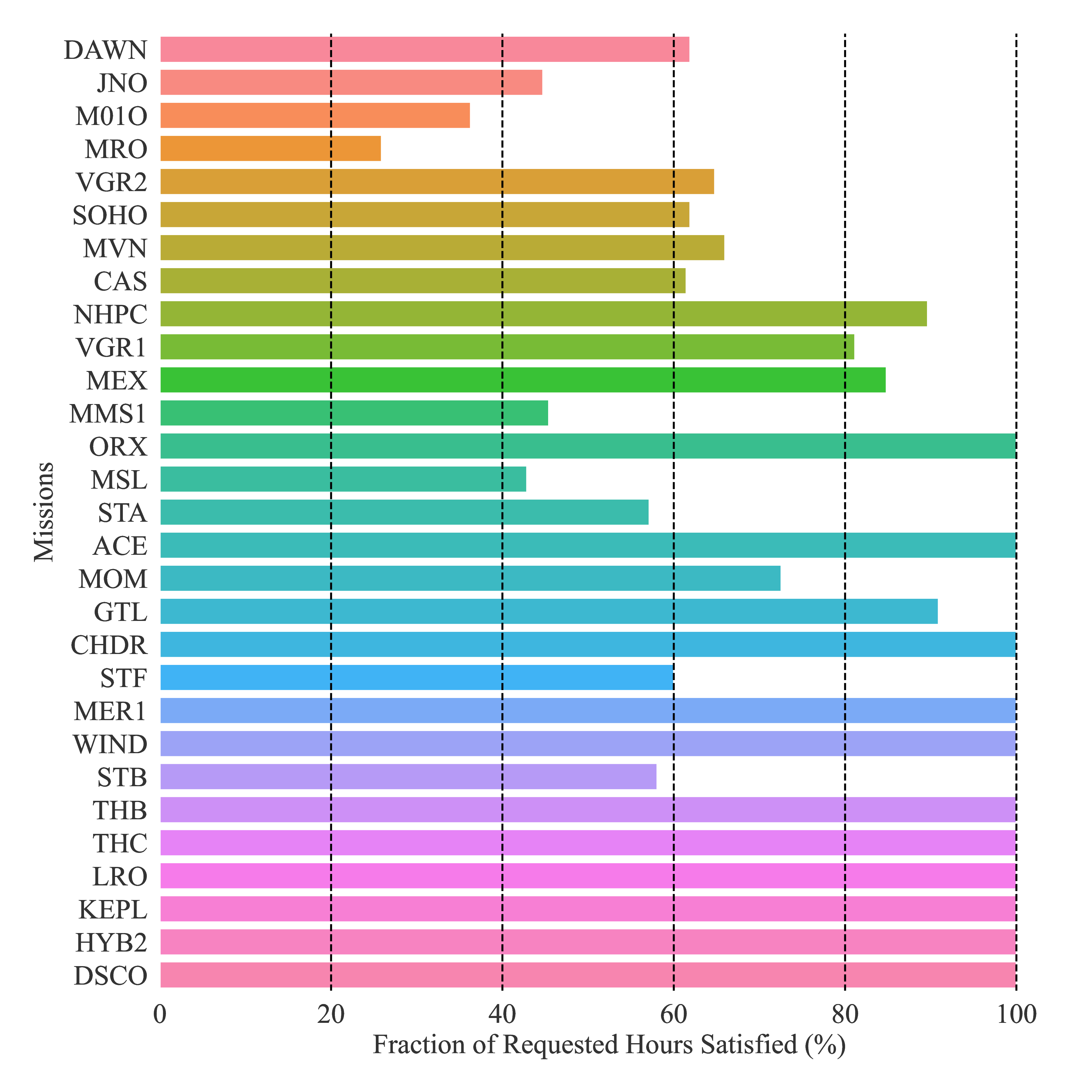}
    \caption{\textbf{Satisfactions for W44 2016 where NHPC is being prioritized with a weight of 5. Comparing with Fig. \ref{fig:results_W44_2016_MILP3}, the satisfaction for NHPC went from 50\% to 90\%, and all the other mission satisfactions went from all above 45\% to all above 26\%. From Fig. \ref{fig:heatmap_scheduled_durations_antennas_MILP_W44_2016}, since NHPC could only be scheduled on 70-m antennas such as DSS-43 and DSS-63, prioritizing it would reduce the amount of scheduled hours for the other missions that need to use the same antennas such as JNO, MRO, or STF.}}
    \label{fig:priority}
\end{figure}

Fig. \ref{fig:priority} shows the results for W44 2016 where NHPC is being prioritized with a weight of 5. Comparing with Fig. \ref{fig:results_W44_2016_MILP3}, the satisfaction for NHPC went from 50\% to 90\%, and all the other mission satisfactions went from all above 45\% to all above 26\%. From Fig. \ref{fig:heatmap_scheduled_durations_antennas_MILP_W44_2016}, since NHPC could only be scheduled on 70-m antennas such as DSS-43 and DSS-63, prioritizing it would reduce the amount of scheduled hours for the other missions that need to use the same antennas such as JNO\footnote{Juno mission  was developed to reveal the story of Jupiter's formation and evolution.}, MRO\footnote{Mars Reconnaissance Orbiter is a spacecraft that was designed to study the geology and climate of Mars, provide reconnaissance of future landing sites, and relay data from surface missions back to Earth.}, or STF\footnote{Splitzer Space Telescope was designed to study the early universe in infrared light.}.

\section{Conclusions and Future Work}\label{sec:conclusion}

In this paper, the Deep Space Network scheduling terminology was first introduced. The problem formulation was then derived with a new set of constraints and this paper presented a new algorithm to improve the user satisfaction.  New results were showed and compared with the previous iterations from \cite{sabol2021deep}, and $\Delta$-MILP outperformed them in terms of satisfaction and feasibility of its schedules.

Parallel developments on DSN scheduling at the Jet Propulsion Laboratory use Deep Reinforcement Learning \cite{goh2021scheduling}. Future work on this research could be to hybridize both approaches. Because is becomes harder to improve the objective function as the schedules get filled, an idea could be to train an agent to make the hardest decisions in a short amount of time. Training this agent could be done using scheduling examples that $\Delta$-MILP would provide by performing longer and more intense searches for this purpose. This hybridization could then help providing better schedules and at a faster pace. Also, comparisons of the $\Delta$-MILP results with other methodologies such as Deep Reinforcement Learning and Quantum Computing solutions will be studied in future publications.  $\Delta$-MILP will be extended to many other applications in the future such as stabilization of balloons \cite{clalimo}, multi-agent motion planning \cite{yun-2020a}, automated data accountability for Mars missions \cite{alimo-2021a} and derivative-free optimization \cite{alimo-2020a}.

\section*{Acknowledgment}
The authors gratefully acknowledge funding from Jet Propulsion Laboratory, California Institute of Technology, under a contract with the National Aeronautics and Space Administration (NASA) in support of this work. The authors thank Alexandre Guillaume, Hamsa Shwetha Venkataram, Chris Mattmann, Joe Lazio, Becky Castano, Thomas Bewley, Alex Sabol, Farhad Kamangar, Dave Hanks, Fred Hadaegh for fruitful discussions and support. 








\newpage
\appendix

\section*{Appendix}
In this appendix, we provide more simulations results in an illustrative manner. The results were reported in Table \ref{tab:comparision_weeks}, but in the following, we only provide the results in a form of Gantt charts for our newly proposed framework $\Delta$-MILP.

\begin{figure*}[ht!]
    
  \includegraphics[width=\textwidth]{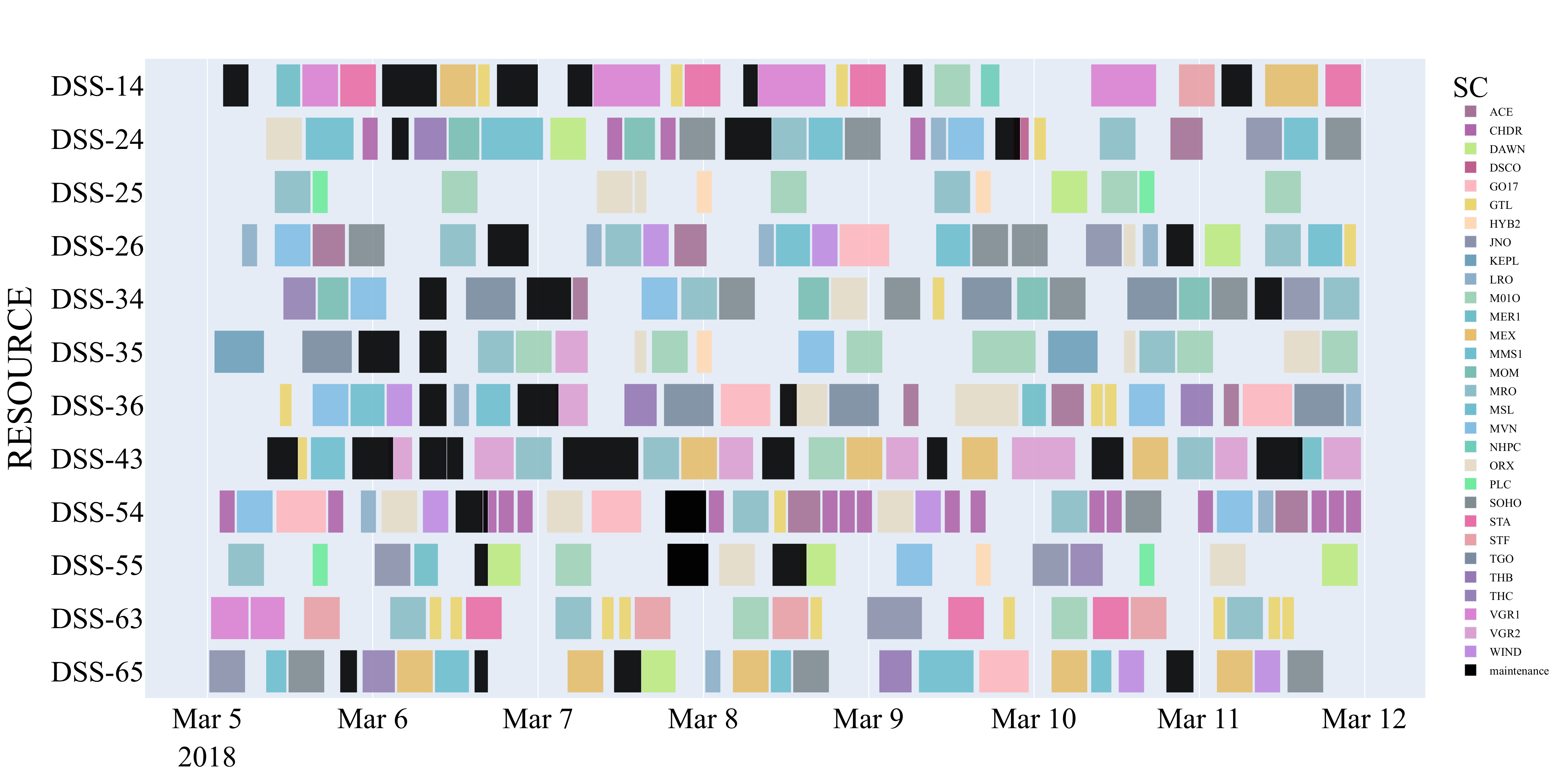}
    \caption{\textbf{Schedule using $\Delta$-MILP for W10 2018.}}
    \label{fig:schedule_DELTA-MILP_W10_2018}
\end{figure*}

\begin{figure*}[ht!]
    
  \includegraphics[width=\textwidth]{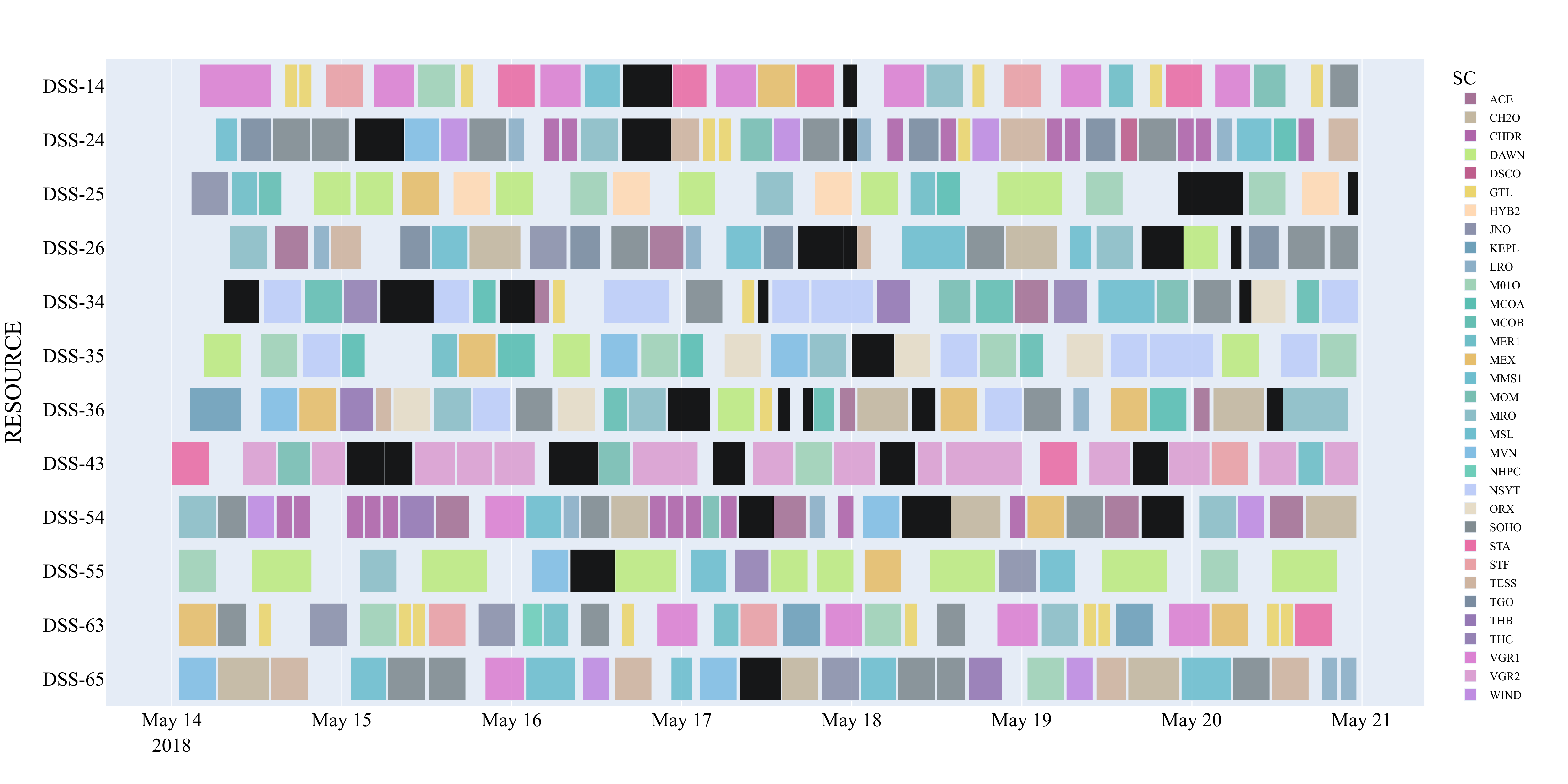}    \caption{\textbf{Schedule using $\Delta$-MILP for W20 2018.}}
    \label{fig:schedule_DELTA-MILP_W20_2018}
\end{figure*}

\begin{figure*}[ht!]
    
  \includegraphics[width=\textwidth]{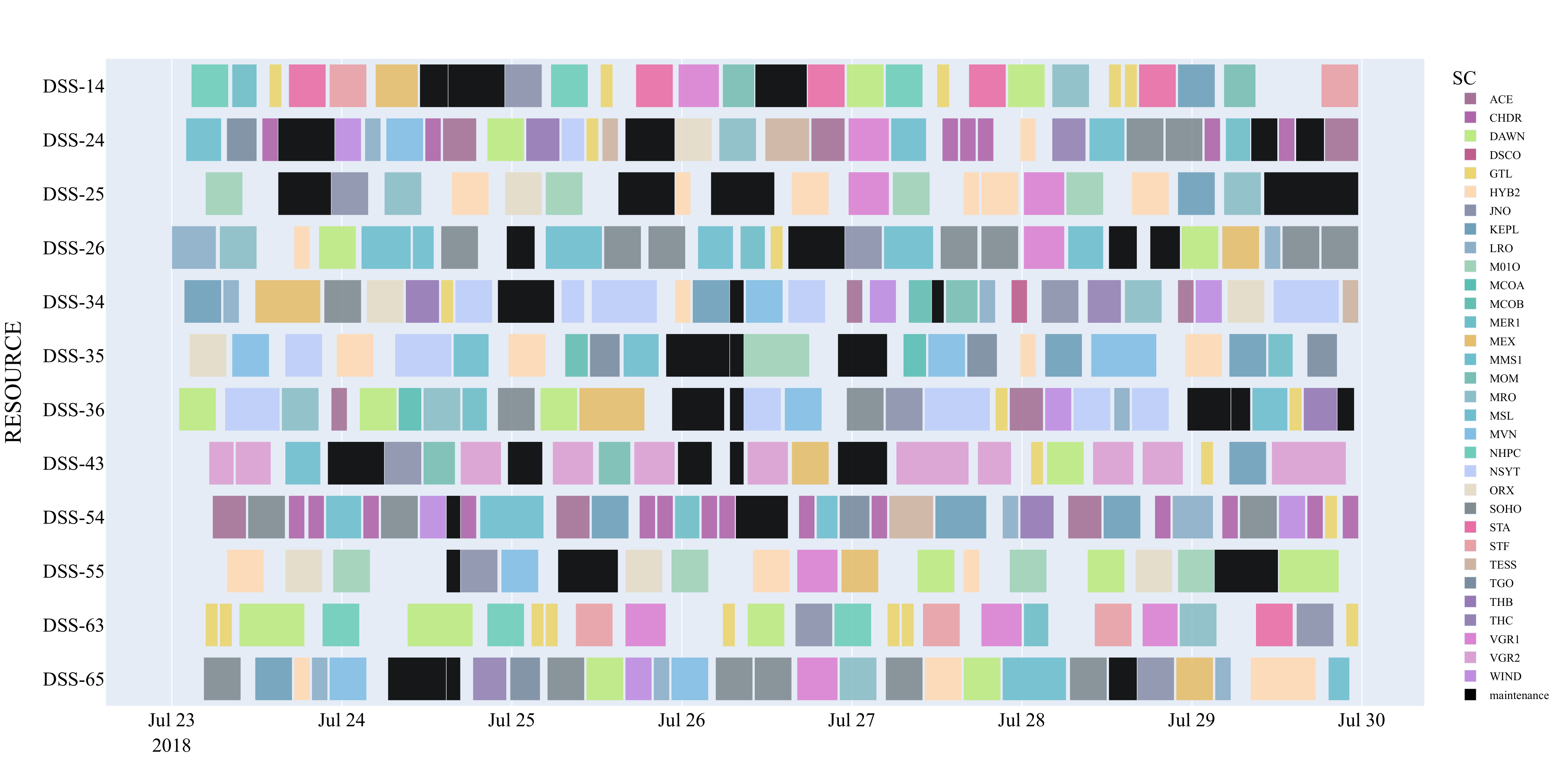}
    \caption{\textbf{Schedule using $\Delta$-MILP for W30 2018.}}
    \label{fig:schedule_DELTA-MILP_W30_2018}
\end{figure*}

\begin{figure*}[ht!]
    
  \includegraphics[width=\textwidth]{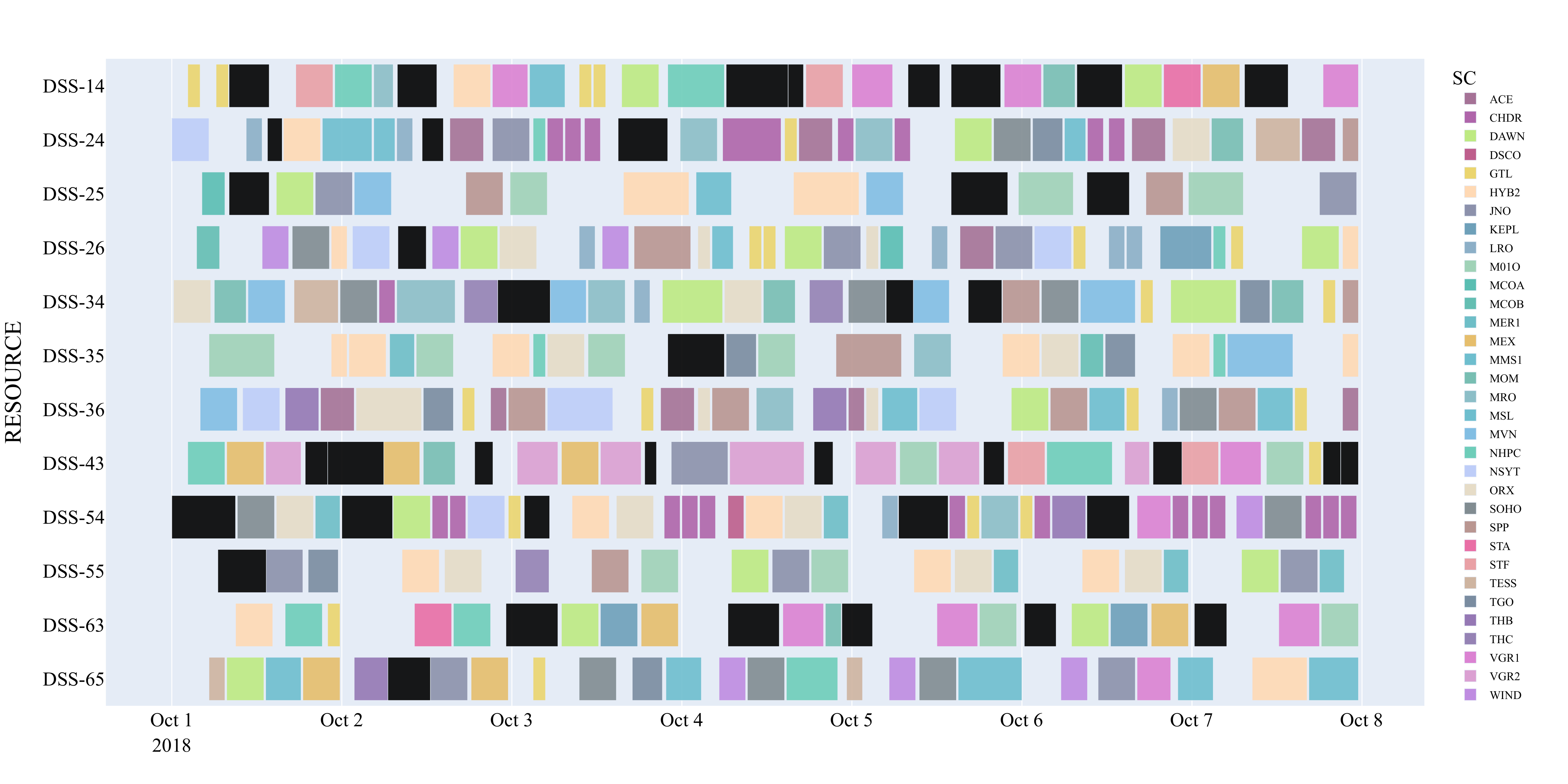}
    \caption{\textbf{Schedule using $\Delta$-MILP for W40 2018.}}
    \label{fig:schedule_DELTA-MILP_W40_2018}
\end{figure*}

\begin{figure*}[ht!]
    
  \includegraphics[width=\textwidth]{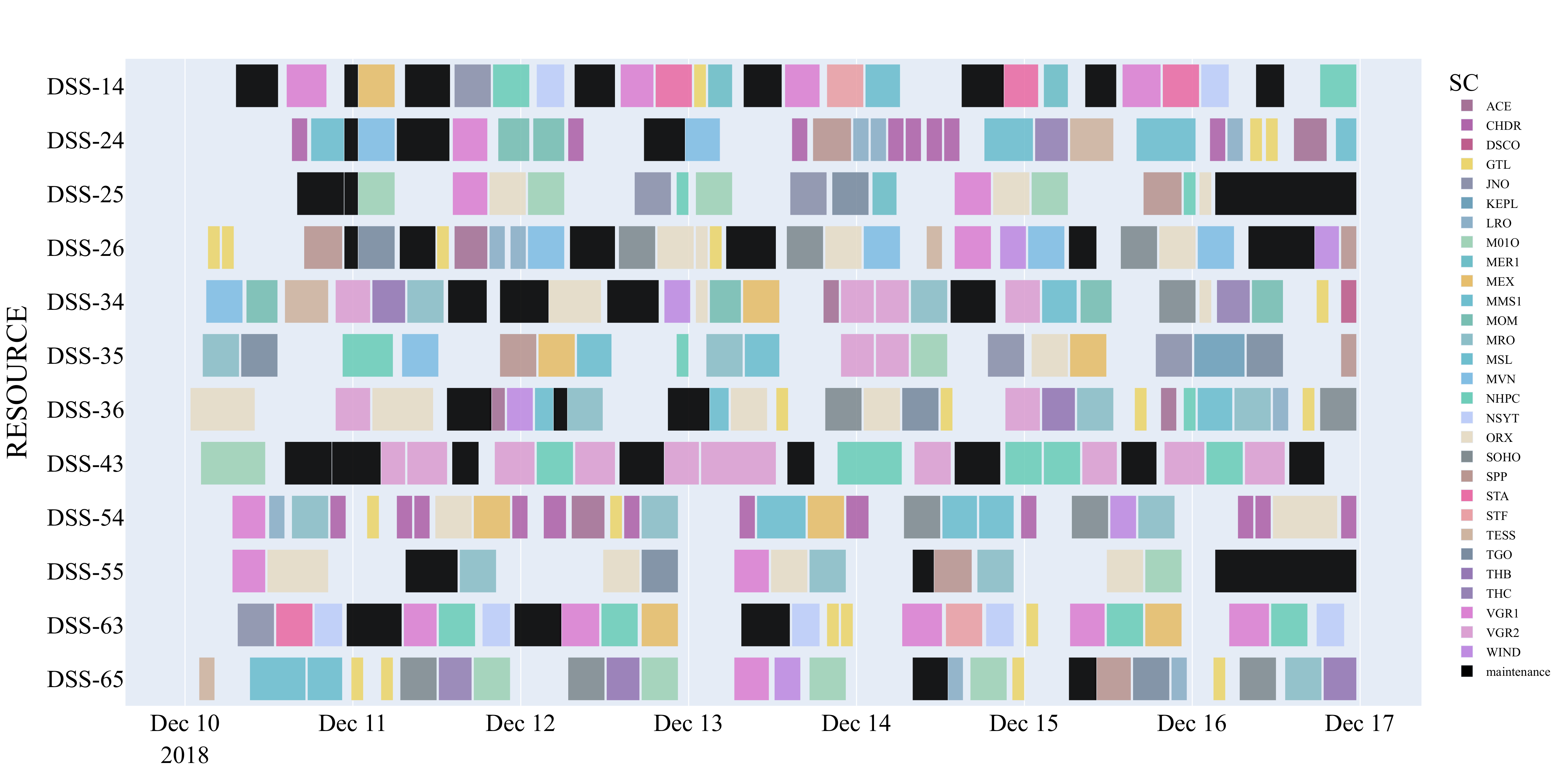}
    \caption{\textbf{Schedule using $\Delta$-MILP for W50 2018.}}
    \label{fig:schedule_DELTA-MILP_W50_2018}
\end{figure*}

\appendix

\end{document}